%% file: main.tex
\documentclass[onefignum,onetabnum]{siamart190516}

\usepackage{hyperref}
\usepackage{fullpage,amsmath,amssymb,enumerate,graphicx,pifont,color,tikz,changepage,yfonts,scalerel,algorithm}
\usepackage{bookmark}

\usetikzlibrary{positioning}
\usetikzlibrary{cd}

\usepackage{booktabs}
\usepackage{multirow}
\usepackage{tabularx}
\usepackage{array}
\newcolumntype{M}[1]{>{\centering\arraybackslash}m{#1}}

\include{tensor_header}

\graphicspath{{fig/}} 

\crefname{section}{Section}{Sections}
\crefname{subsection}{Section}{Sections}

\title{CP Decomposition for Tensors via Alternating Least Squares with QR Decomposition}
\author{Rachel Minster, Irina Viviano, Xiaotian Liu, Grey Ballard}

\begin{document}

\maketitle

\begin{abstract}
\addcontentsline{toc}{section}{Abstract}
The CP tensor decomposition is used in applications such as machine learning and signal processing to discover latent low-rank structure in multidimensional data.
Computing a CP decomposition via an alternating least squares (ALS) method reduces the problem to several linear least squares problems.
The standard way to solve these linear least squares subproblems is to use the normal equations, which inherit special tensor structure that can be exploited for computational efficiency.
However, the normal equations are sensitive to numerical ill-conditioning, which can compromise the results of the decomposition. 
In this paper, we develop versions of the CP-ALS algorithm using the QR decomposition and the singular value decomposition (SVD), which are more numerically stable than the normal equations, to solve the linear least squares problems. 
Our algorithms utilize the tensor structure of the CP-ALS subproblems efficiently, have the same complexity as the standard CP-ALS algorithm when the rank is small, and are shown via examples to produce more stable results when ill-conditioning is present. 
Our MATLAB implementation achieves the same running time as the standard algorithm for small ranks, and we show that the new methods can obtain lower approximation error and more reliably recover low-rank signals from data with known ground truth.
\end{abstract}


%

\section{Introduction}
\label{sec:intro}

The CANDECOMP/PARAFAC or canonical polyadic (CP) decomposition for multidimensional data, or tensors, is a popular tool for analyzing and interpreting latent patterns that may be present in multidimensional data. 
One of the most popular methods used to compute a CP decomposition is the alternating least squares (CP-ALS) approach, which solves a series of linear least squares problems \cite{tensor-toolbox,tensorlab,Tensorly}. 
To solve these linear least squares problems, CP-ALS uses the normal equations, which are well known to be sensitive to roundoff error for moderately ill-conditioned matrices. 
We propose to use a more stable approach, where we solve the linear least squares problems using the QR decomposition instead. 

Consider the standard linear least squares problem with multiple outputs, i.e.,
$\min_{\M{X}} \| \M{AX}-\M{B}\|_F.$
The normal equations for this problem are $\M{A}^\top \M{A}\M{X} = \M{A}^\top \M{B}$. 
Given the compact/thin QR decomposition $\M{A}=\M{Q}\M{R}$, the more numerically stable solution is computed from $\M{R}\M{X} = \M{Q}^\top \M{B}$.
When $\M{A}$ is tall and skinny and $\M{B}$ has few columns, the normal equations approach has about half the cost of the QR-based approach because the dominant costs are computing $\M{A}^\top \M{A}$ and the QR decomposition, respectively.
However, if $\M{B}$ has many more columns than $\M{A}$, then the dominant costs are those of computing $\M{A}^\top \M{B}$ and $\M{Q}^\top\M{B}$, which are equivalent when $\M{Q}$ is formed explicitly.
In this case, the QR-based approach is more numerically stable and requires practically no more computation compared to the normal equations approach.
Furthermore, for rank-deficient problems, the QR-based approach can be cheaply extended to use the SVD to solve the least squares problem, computing the minimum norm solution among the set of solutions with equivalent residual norm.

As we describe in more detail in \cref{sec:background}, when solving linear least squares problems within CP-ALS, $\M{A}$ corresponds to a Khatri-Rao product of factor matrices, and $\M{B}$ corresponds to the transpose of a matricized tensor. 
In this case, the number of columns of $\M{A}$ is the rank of the CP decomposition, and the number of columns of $\M{B}$ corresponds to one of the tensor dimensions.
The normal equations are particularly convenient within CP-ALS because the Khatri-Rao structure of $\M{A}$ can be exploited to compute $\M{A}^\top\M{A}$ very efficiently, and thus, even for large ranks, the dominant cost is computing $\M{A}^\top \M{B}$, whose transpose is known as the matricized-tensor times Khatri-Rao product (MTTKRP).
The MTTKRP is a well-studied and well-optimized computation because of its importance for the performance of CP-ALS and other gradient-based optimization algorithms for CP \cite{BR20,NL+19,PTC13a,splat}.

In \cref{ssec:qr-solve}, we present a QR-based approach to solving the linear least squares problems within CP-ALS.
In order to achieve comparable computational complexity with the normal equations approach, we exploit the Khatri-Rao structure of $\M{A}$ in the computation of the QR decomposition as well as $\M{Q}^\top\M{B}$.
In particular, we show that the QR decomposition of a Khatri-Rao product of matrices can be computed efficiently from the QR decompositions of the individual factor matrices.
Using the structure of the orthonormal component $\M{Q}$, the computation of $\M{Q}^\top\M{B}$ involves multiple tensor-times-matrix (Multi-TTM) products.
We prove in \cref{ssec:costs} that when the rank is small relative to the tensor dimensions, the Multi-TTM is the dominant cost, and it has the same leading-order complexity as MTTKRP for dense tensors.
Multi-TTM is also a well-optimized tensor computation, as it is important for the computation of Tucker decompositions \cite{BKK20,CC+17,KU16,smith-tucker}.

When the rank is comparable to or much larger than the tensor dimensions, the QR-based approach can require significantly more computation time than standard CP-ALS using the normal equations.
Both the QR decomposition and application of the orthonormal component have costs that are lower order only when the rank is small.
In this case, the numerical stability provided by QR comes at the expense of performance.

We demonstrate the performance and accuracy of our methods using several example input tensors in \cref{sec:examples}.
Our MATLAB implementation of the algorithms uses the Tensor Toolbox \cite{tensor-toolbox}, and we compare against the CP-ALS algorithm implemented in that library.
In \cref{ssec:perf}, we validate the theoretical complexity analysis and show that there is no increase in per-iteration costs when the rank is small, and we demonstrate with a time breakdown which of the computations become bottlenecks as the rank grows larger relative to the tensor dimensions.
To illustrate the differences in accuracy, we present two sets of examples in \cref{ssec:collinear,ssec:sine} that lead to ill-conditioned subproblems and show that the instability of the normal equations can lead to degradation of approximation error, recovery of low-rank signal when the ground truth is known, and in some cases, convergence of the overall algorithm.

We conclude in \cref{sec:conclusion} that using the QR-based approaches to solve the CP-ALS subproblems increases the robustness of the overall algorithm without sacrificing performance in the typical case of small ranks.
However, due to the complexity of the algorithm and the extra computation that becomes significant for large ranks, we envision a CP-ALS solver that uses the fast-and-inaccurate normal equations approach by default and falls back on the accurate-but-possibly-slow SVD approach when necessary.
For problems that do not involve ill-conditioning, which is the case for many tensors representing noisy data, the normal equations are sufficient for obtaining accurate solutions.
However, when ill-conditioning degrades the accuracy of solutions computed from the normal equations, we show that it is possible to obtain the stability of the SVD with feasible computational cost.

\section{Background}\label{sec:background}
In this section, we first review typical methods for solving linear least squares problems. We also discuss the relevant information regarding tensors and the CP decomposition, focusing on the alternating least squares (CP-ALS) algorithm. We also briefly describe an optimization approach for CP using Gauss-Newton. 

\subsection{Linear Least Squares Methods}
A common approach to solving linear least squares problems is by solving the associated normal equations.
When applied to ill-conditioned problems, however, using the normal equations results in numerical instability. 
More numerically stable methods to solve least squares problems include using the QR decomposition or the SVD. 

Consider a least squares problem of the form 
\begin{equation*}
	\min_{\M{X}} \| \M{B} - \M{X} \M{A}^\top \|_F.
\end{equation*}
Note that the coefficient matrix appears to the right of the variable matrix rather than the left (as appears in \cref{sec:intro}) in order to match the form of the CP-ALS subproblems described below.
The normal equations for this problem are $\M{X} \M{A}^\top\M{A} = \M{B}\M{A}$, which is equivalent to $\M{A}^\top\M{A} \M{X}^\top = \M{A}^\top\M{B}^\top$.
To solve this least squares problem using QR, we first compute the compact/thin QR factorization of $\M{A} = \M{QR}$, so that $\M{Q}$ has the same dimensions as $\M{A}$ and $\M{R}$ is square. We then apply $\M{Q}$ to matrix $\M{B}$ on the right, and use the result to solve the triangular system $\M{XR}^\top = \M{B} \M{Q}$ for $\M{X}$.
When the coefficient matrix $\M{A}$ is tall and skinny, the QR approach can be cheaply extended to use the SVD. 
Given the QR factorization of $\M{A}$, we compute the SVD of $\M{R} = \M{U \Sigma V}^\top$. 
If we apply $\M{U}$ to $\M{B}\M{Q}$ on the right, we can then solve the system $\M{Y}\M{\Sigma}  = \M{B}\M{Q}\M{U}$ for $\M{Y}$ and compute $\M{X} =  \M{Y} \M{V}^\top$.
If $\M{A}$ is numerically low rank, we can solve the system using the pseudoinverse of $\M{\Sigma}$ to find the minimum-norm solution to the original problem.
For more details on methods to solve linear least squares problems, see \cite{Demmel97,GVL13,trefethen-book}.

\subsection{Tensor Notation and Preliminaries}
Throughout this paper, we follow the notation from \cite{intro-paper}. A scalar is denoted by lowercase letters, e.g. $a$, while vectors are denoted by boldface lowercase letters, e.g. $\V{a}$. Matrices are denoted by boldface uppercase letters, e.g. $\M{A}$, and tensors are denoted by boldface uppercase calligraphic letters, e.g. $\T{X}$. We use the \textsf{MATLAB} notation $\M{A}(i,:)$ to refer to the $i^{\text{th}}$ row of $\M{A}$ and $\M{A}(:,j)$ to refer to the $j^{\text{th}}$ column of $\M{A}$. 

\paragraph{Matrix Products}
We define three matrix products that will appear in our algorithms. First, the Kronecker product of two matrices $\M{A} \in \mb{R}^{I \times J}$ and $\M{B} \in \mb{R}^{K \times M}$ is denoted $\M{A} \Kron \M{B} \in \mb{R}^{(IK) \times (JM)}$ matrix with entries $[\M{A} \Kron \M{B}](K(i-1)+k, M(j-1) + \ell) = \M{A}(i,j)\M{B}(k,\ell)$. 

The Khatri-Rao product of matrices $\M{A} \in \mb{R}^{I \times K}$ and $\M{B} \in \mb{R}^{J \times K}$ is denoted $\M{A} \Khat \M{B} \in \mb{R}^{(IJ) \times K}$ matrix with columns $[\M{A} \Khat \M{B}](:,k) = \M{A}(:,k) \Kron \M{B}(:,k)$ for every $k = 1,\dots, K$.
We present a property regarding the Khatri-Rao product of a product of two matrices, which will be useful in a later section. 
Let $\M{A} \in \mb{R}^{K \times J}, \M{B} \in \mb{R}^{I \times J}, \M{C} \in \mb{R}^{I \times K}$, and $\M{D} \in \mb{R}^{J \times I}$ be four matrices. Then,
\begin{equation}\label{eq:kr-id}
(\M{C} \Kron \M{D}) (\M{A} \Khat \M{B}) = (\M{C} \M{A}) \Khat (\M{D} \M{B}).
\end{equation}
Finally, the Hadamard product of two matrices $\M{A} \in \mb{R}^{I \times J}$ and $\M{B} \in \mb{R}^{I \times J}$ is the elementwise product denoted as $\M{A} * \M{B} \in \mb{R}^{I \times J}$, with entries $[\M{A} * \M{B}](i,j) = \M{A}(i,j)\M{B}(i,j)$.

\paragraph{Tensor Components}
As generalizations to matrix rows and columns, tensors have mode-$j$ fibers, which are vectors formed by fixing all but one index of a tensor. Similarly, slices of a tensor are two-dimensional sections, formed by fixing all but two indices of a tensor. 

\paragraph{Tensor Operations}
There are two major tensor operations we will frequently use throughout this paper, namely matricization and multiplying a tensor with a matrix.
The $n$-mode matricization, or unfolding, of a tensor $\T{X} \in  \mathbb{R}^{I_1 \times \cdots \times I_N}$, denoted $\Mz{X}{n} \in \mb{R}^{I_n \times \left(\prod_{j \neq n} I_j \right)}$, and the matrix $\Mz{X}{n}$ is formed so that the columns are the mode-$n$ fibers of $\T{X}$.

Also useful is the tensor-times-matrix (TTM) multiplication, or the $n$-mode product of a tensor and matrix. Let $\T{X} \in \mathbb{R}^{I_1 \times \cdots \times I_N}$ and $\M{U} \in \mathbb{R}^{J \times I_n}$. The resulting tensor $\T{Y} = \T{X} \times_n \M{U} \in \mb{R}^{I_1 \times \cdots \times I_{n-1} \times J \times I_{n+1} \times \cdots \times I_N}$, and it has entries $ \T{Y}(i_1,\dots,i_{n-1},j,i_{n+1},\dots,i_N) = \sum_{i_n = 1}^{I_n} \T{X}(i_1, \dots, i_N) \M{U}(j, i_n)$.
The TTM can also be computed via the mode-$n$ matricization $\Mz{Y}{n} = \M{U} \Mz{X}{n}$. 

Multiplying an $N$-mode tensor by multiple matrices in distinct modes is known as Multi-TTM.
The computation can be performed using a sequence of individual mode TTMs, and they can be done in any order.
In particular, we will be interested in the case where we multiply an $N$-mode tensor by matrices $\M{U}_j$ in every mode except $n$, denoted $\T{Y} = \T{X} \times_1 \M{U}_1 \dots \times_{n+1} \M{U}_{n-1} \times_{n+1} \M{U}_{n+1} \dots \times_N \M{U}_N$. 
This is expressed in the mode-$n$ matricization as
\begin{equation}
\label{eq:multittm}
	\Mz{Y}{n} = \Mz{X}{n}(\M{U}_N \otimes \dots \otimes \M{U}_{n+1} \otimes \M{U}_{n-1} \otimes \dots \otimes \M{U}_1)^\top.
\end{equation}

\subsection{CP-ALS Algorithm}
We now detail both the CP decomposition as well as the alternating least squares approach (CP-ALS), one of the most popular algorithms used to compute the CP decomposition.

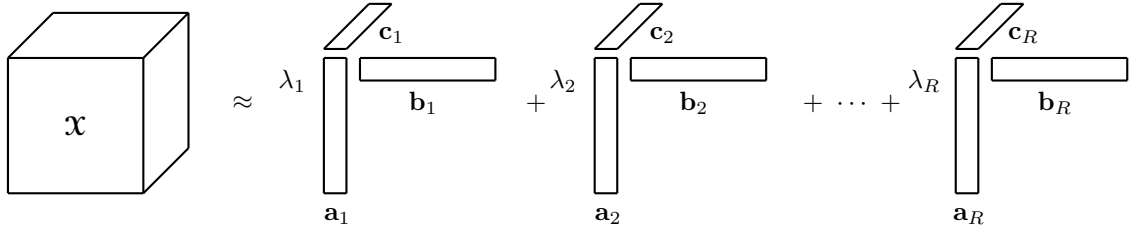
\begin{figure}[H]
\centering
\begin{tikzpicture}[scale=0.6]

\draw node (X) at (1.5, 2.5) {$\T{X}$};
\draw[thick] (0,1) -- (3,1);
\draw[thick] (0,1) -- (0,4);
\draw[thick] (0,4) -- (3,4);
\draw[thick] (3,4) -- (3,1);
\draw[thick] (0,4) -- (1,5);
\draw[thick] (1,5) -- (4,5);
\draw[thick] (3,4) -- (4,5);
\draw[thick] (4,5) -- (4,2);
\draw[thick] (4,2) -- (3,1);

\draw node at (5.2, 3) {$\approx$};

\draw node at (6.3, 3.5) {$\lambda_1$};

\draw node at (7.3, 0.5) {$\V{a}_1$};
\draw[thick] (7,1) -- (7.5,1);
\draw[thick] (7,1) -- (7,4);
\draw[thick] (7,4) -- (7.5,4);
\draw[thick] (7.5,1) -- (7.5,4);

\draw node at (9.2, 3) {$\V{b}_1$};
\draw[thick] (7.8,3.5) -- (7.8,4);
\draw[thick] (7.8,3.5) -- (10.8,3.5);
\draw[thick] (7.8,4) -- (10.8,4);
\draw[thick] (10.8,4) -- (10.8,3.5);

\draw node at (8.5, 4.5) {$\V{c}_1$};
\draw[thick] (7,4.2) -- (7.5,4.2);
\draw[thick] (7,4.2) -- (8,5.2);
\draw[thick] (8,5.2) -- (8.5,5.2);
\draw[thick] (7.5,4.2) -- (8.5,5.2);

\draw node at (11.7, 3) {$+$};

\draw node at (12.3, 3.5) {$\lambda_2$};

\draw node at (13.3, 0.5) {$\V{a}_2$};
\draw[thick] (13,1) -- (13.5,1);
\draw[thick] (13,1) -- (13,4);
\draw[thick] (13,4) -- (13.5,4);
\draw[thick] (13.5,1) -- (13.5,4);

\draw node at (15.2, 3) {$\V{b}_2$};
\draw[thick] (13.8,3.5) -- (13.8,4);
\draw[thick] (13.8,3.5) -- (16.8,3.5);
\draw[thick] (13.8,4) -- (16.8,4);
\draw[thick] (16.8,4) -- (16.8,3.5);

\draw node at (14.5, 4.5) {$\V{c}_2$};
\draw[thick] (13,4.2) -- (13.5,4.2);
\draw[thick] (13,4.2) -- (14,5.2);
\draw[thick] (14,5.2) -- (14.5,5.2);
\draw[thick] (13.5,4.2) -- (14.5,5.2);

\draw node at (18.7, 3) {$+ \ \cdots \ +$};

\draw node at (20.3, 3.5) {$\lambda_R$};

\draw node at (21.3, 0.5) {$\V{a}_R$};
\draw[thick] (21,1) -- (21.5,1);
\draw[thick] (21,1) -- (21,4);
\draw[thick] (21,4) -- (21.5,4);
\draw[thick] (21.5,1) -- (21.5,4);

\draw node at (23.2, 3) {$\V{b}_R$};
\draw[thick] (21.8,3.5) -- (21.8,4);
\draw[thick] (21.8,3.5) -- (24.8,3.5);
\draw[thick] (21.8,4) -- (24.8,4);
\draw[thick] (24.8,4) -- (24.8,3.5);

\draw node at (22.5, 4.5) {$\V{c}_R$};
\draw[thick] (21,4.2) -- (21.5,4.2);
\draw[thick] (21,4.2) -- (22,5.2);
\draw[thick] (22,5.2) -- (22.5,5.2);
\draw[thick] (21.5,4.2) -- (22.5,5.2);
\end{tikzpicture}
\caption{CP decomposition of rank $R$ for a three-dimensional tensor $\T{X}$ \label{fig:3d-cp-decomp}}
\end{figure}

\paragraph{CP Decomposition}
The aim of the CP decomposition is to represent a tensor as a sum of rank-one components. For an $N$-mode tensor $\T{X} \in \mb{R}^{I_1 \times \dots \times I_N}$, the rank-$R$ CP decomposition of $\T{X}$ is the approximation
\begin{equation}\label{eq:cp}
	\hat{\T{X}} = \sum_{r=1}^R \lambda_r \, \V{a}_r^{(1)} \circ \V{a}_r^{(2)} \circ \dots \circ \V{a}_r^{(N)}.
\end{equation}
where $\V{a}_r^{(n)} \in \mb{R}^{I_n}$ are unit vectors with weight vector $\V{\lambda} \in \mb{R}^{R}$, and $\circ$ denotes the outer product. The collection of all $\V{a}_r^{(n)}$ vectors for each mode is called a factor matrix, i.e., $\M{A}_{n} = \bmat{ \V{a}_1^{(n)} & \V{a}_2^{(n)} & \dots & \V{a}_r^{(n)}} \in \mb{R}^{I_n \times R}.$
A visualization of the three-dimensional version of this representation is in \cref{fig:3d-cp-decomp}.

Using the notation $\M{\hat{A}}_n = \M{A}_n \cdot \text{diag}(\M{\lambda})$, we can express the mode-$n$ matricization of $\hat{\T{X}}$  as
\begin{equation}\label{eq:cp_mat}
	\Mz{\hat{X}}{n} = \M{\hat{A}}_{n}(\M{A}_{N} \odot \dots \odot \M{A}_{n+1} \odot \M{A}_{n-1} \odot \dots \odot \M{A}_{1})^\top = \M{A}_{n}\M{Z}_n^\top,
\end{equation}
letting $\M{Z}_{n} = \M{A}_{N} \odot \dots \odot \M{A}_{n+1} \odot \M{A}_{n-1} \odot \dots \odot \M{A}_{1}$.

\paragraph{CP-ALS Algorithm}
In order to compute the CP-decomposition, the CP-ALS algorithm solves a least squares problem for the matricized tensors $\M{X}_{(n)}$ and $\hat{\M{X}}_{(n)}$ along each mode. 
For mode $n$, we fix every factor matrix except $\M{\hat{A}}_{n}$ and then solve for $\Mn{\hat{A}}{n}$. 
This process is repeated by alternating between modes until some termination criteria is met. 
We solve the linear least squares problem 
\begin{equation}\label{eq:ls_n}
	\min_{\M{\hat{A}}_{n}} \| \Mz{X}{n} - \M{\hat{A}}_{n}\M{Z}_{n}^\top \|_F,
\end{equation}
using the representation in \cref{eq:cp_mat}. 
The linear least squares problem from \cref{eq:ls_n} is typically solved using the normal equations, i.e.,
\begin{equation*}
	\Mz{X}{n} \M{Z}_{n} = \M{\hat{A}}_{n} ( \M{Z}_{n}^\top \M{Z}_{n}).
\end{equation*}
The coefficient matrix $\M{Z}_{n}^\top \M{Z}_{n}$ is computed efficiently as $\M{A}_1^\top \M{A}_1 \Hada \dots \Hada \M{A}_{n-1}^\top \M{A}_{n-1} \Hada \M{A}_{n+1}^\top \M{A}_{n+1} \Hada \dots \Hada \M{A}_{N}^\top \M{A}_{N}$.
We obtain the desired factor matrix $\M{A}_n$ by normalizing the columns of $\M{\hat{A}}_n$ and updating the weight vector $\M{\lambda}$. This approach is detailed in \cref{alg:cp-als}.

\begin{algorithm}[!ht]
  \caption{CP-ALS}
  \label{alg:cp-als}
  \begin{algorithmic}[1]\footnotesize
    \Function{$[\bm{\lambda},\set{\M{A}_{n}}]=$ CP-ALS}{$\T{X},R$}\Comment{$\T{X}\in\mathbb{R}^{I_1\times \cdots \times I_N}$}
    \State \label{line:cpals:init}Initialize factor matrices $\M{A}_{1}, \dots, \M{A}_{N}$
    \State \label{line:gram}Compute Gram matrices $\Mn{G}{1} = \M{A}_1^\top \Mn{A}{1}, \ldots, \Mn{G}{N} = \M{A}_N^\top \Mn{A}{N}$
    \Repeat
    \For{$n=1,\dots, N$}
      \State $\M{S}_n \gets \Mn{G}{N} \Hada \cdots \Hada \Mn{G}{n+1} \Hada \Mn{G}{n-1} \Hada \cdots \Hada \Mn{G}{1}$\label{line:cpals:Gram}
      \State \label{line:cpals:KR}$\M{Z}_{n} \gets \M{A}_{N}\odot \cdots  \odot \M{A}_{n+1}\odot \M{A}_{n-1} \odot \cdots \odot \M{A}_{1}$
      \State \label{line:mttkrp}$\M{M}_n \gets \M{X}_{(n)}\M{Z}_{n}$\Comment{MTTKRP}
      \State \label{line:cpals:solve}Solve $\M{\hat{A}}_{n}\M{S}_n = \M{M}_n$ for $\M{\hat{A}}_{n}$ via Cholesky (CP-ALS) or SVD (CP-ALS-PINV) \Comment{Normal equations}
      \State Normalize columns of $\M[\hat]{A}_{n}$ to obtain $\M{A}_n$ and $\V{\lambda}$
      \State Recompute Gram matrix $\Mn{G}{n} = \M{A}_{n}^\top \Mn{A}{n}$ for updated factor matrix $\Mn{A}{n}$
    \EndFor
    \Until termination criteria met
    \State \textbf{return} $\bm{\lambda}$, factor matrices $\set{\M{A}_{n}}$
    \EndFunction
  \end{algorithmic}
\end{algorithm}

Note that in \cref{line:cpals:solve}, the standard approach (implemented in the Tensor Toolbox \cite{tensor-toolbox}) is to solve $\M{A}_{n}\M{S} = \M{M}$ for $\M{A}_{n}$ using a Cholesky decomposition, which is implemented in MATLAB using the backslash operator. 
An alternative method for solving the linear system is to use the SVD, which is implemented in MATLAB using the \texttt{pinv} function. 
We call the algorithm obtained by taking this alternate approach CP-ALS-PINV.

\subsection{Gauss-Newton Optimization Approach}\label{ssec:gn}
An alternate approach to alternating least squares is to minimize the CP model ``all at once'' using optimization techniques. 
Instead of the linear least squares problem from \cref{eq:ls_n}, this approach solves the nonlinear least squares problem $\min \| \T{X} - \hat{\T{X}} \|_F^2$, subject to factor matrices $\M{A}_j$ with $j = 1,\dots, N$, where $\hat{\T{X}}$ is defined in \cref{eq:cp}. 
Gauss-Newton attempts to minimize this nonlinear residual by using a linear Taylor series approximation at each iteration and minimizing that function using standard linear least squares methods. 
Typically, these linear least squares problems are solved via the normal equations as follows.  For residual $\V{r} = \text{vec}{(\T{X} - \hat{\T{X}})}$ and $\M{J}$ the Jacobian of $\V{r}$, the normal equations to be solved are $\M{J}^\top \M{J} = \M{J} \V{r}$. While CP-ALS is typically fast and easy to implement, the Gauss-Newton approach is beneficial as it can converge quadratically if the residual is small. 
We use the implementation of Gauss-Newton in Tensorlab \cite{tensorlab} in the experiments in \cref{sec:examples}.
For more details on this approach, see \cite{singhGN,sorberGN,vervlietGN}.

\section{Proposed Methods}\label{sec:cp-als-qr}

\subsection{CP-ALS-QR Algorithms}\label{ssec:qr-solve} 

In our proposed CP-ALS approach, we incorporate the more stable QR decomposition and SVD and avoid using the normal equations.
Suppose $\T{X} \in \mathbb{R}^{I_1 \times \cdots \times I_N}$ is an $N$-dimensional tensor and we wish to approximate a solution $\T{\hat{X}} = [\![ \V{\lambda}; \Mn{A}{1}, \ldots, \Mn{A}{N} ]\!]$ of rank $R$. 
Recall our linear least-squares problem from \cref{eq:ls_n}.
The key to the efficiency of our algorithm is to form the QR decomposition of $\M{Z}_{n}$ using the Khatri-Rao structure as follows. 
As a first step, we compute compact QR factorizations of each individual factor matrix, so thatLet each factor matrix $\M{A}_j$ have compact QR factorization $\M{A}_j = \M{Q}_j \M{R}_j$. 
Then
\begin{equation*}
\begin{aligned}
	\M{Z}_n &= \M{A}_N \odot \dots \odot \M{A}_{n+1} \odot \M{A}_{n-1} \odot \dots \odot \M{A}_1 \\
	&= \M{Q}_n \M{R}_N \odot \dots \odot \M{Q}_{n+1}\M{R}_{n+1} \odot \M{Q}_{n-1} \M{R}_{n-1} \odot \dots \odot \M{Q}_1 \M{R}_1 \\
	&= (\M{Q}_N \otimes \dots \otimes \M{Q}_{n+1} \otimes \M{Q}_{n-1} \otimes \dots \otimes \M{Q}_1)\underbrace{(\M{R}_N \odot \dots \odot \M{R}_{n+1} \odot \M{R}_{n-1} \odot \dots \odot \M{R}_1)}_{\M{V}_n},
\end{aligned}
\end{equation*}
where the last equality comes from \cref{eq:kr-id}. We then compute the QR factorization of the Khatri-Rao product $\M{V}_n=\M{R}_N \odot \dots \odot \M{R}_{n+1} \odot \M{R}_{n-1} \odot \dots \odot \M{R}_1 = \M{Q}_0 \M{R}_0$. This allows us to express the QR of $\M{Z}_n$ as 
\begin{equation}\label{eq:Z_qr}
\begin{aligned}
	\M{Z}_n &= (\M{Q}_N \otimes \dots \otimes \M{Q}_{n+1} \otimes \M{Q}_{n-1} \otimes \dots \otimes \M{Q}_1)(\M{R}_N \odot \dots \odot \M{R}_{n+1} \odot \M{R}_{n-1} \odot \dots \odot \M{R}_1) \\
	&= \underbrace{(\M{Q}_N \otimes \dots \otimes \M{Q}_{n+1} \otimes \M{Q}_{n-1} \otimes \dots \otimes \M{Q}_1) \M{Q}_0}_{\M{Q}} \underbrace{\M{R}_0}_{\M{R}}.
\end{aligned}
\end{equation}
Note that $\M{Q}$ has orthonormal columns, and $\M{R}$ is upper triangular. 
The Khatri-Rao product of triangular matrices $\M{V}_n$ has sparse structure. 
One column is dense, but the rest have many zeros. 
As discussed in \cref{ssec:sparsity}, this sparse structure can be exploited while implementing the QR of $\M{V}_n$, but our current implementation treats it as dense.

Once the QR is computed, we apply the representation in \cref{eq:Z_qr} to our least squares problem to obtain
\begin{equation}\label{eq:ls-qr-step}
 \min_{\Mn{\hat{A}}{n}} \| \Mz{X}{n} - \Mn{\hat{A}}{n} \M{R}^\top \M{Q}_0^\top (\M{Q}_{N} \otimes \dots \otimes \M{Q}_{n+1} \otimes \M{Q}_{n-1} \otimes \dots \otimes \M{Q}_1 )^\top \|_F,
\end{equation}
with $\M{R} \in \mb{R}^{R \times R}$, $\M{Q}_0 \in \mb{R}^{R^{N-1} \times R}$, and each $\M{Q}_{i} \in \mb{R}^{I_i \times R}$. We are solving for the unnormalized factor matrix $\Mn{\hat{A}}{n} = \Mn{A}{n} \cdot \text{diag}(\M{\lambda})$.

Next, we apply the product $\M{Q}_N \otimes \dots \otimes \M{Q}_{n+1} \otimes \M{Q}_{n-1} \otimes \dots \otimes \M{Q}_1$ to $\Mz{X}{n}$ on the right via the Multi-TTM $\T{Y} = \T{X} \times_1 \M{Q}_1^\top \dots \times_{n-1} \M{Q}_{n-1}^\top \times_{n+1} \M{Q}_{n+1}^\top \dots \times_N \M{Q}_N^\top$, following \cref{eq:multittm}. 
Our least-squares problem then becomes
\begin{equation}\label{eq:ls-qr-ttm-step}
\min_{\Mn{\hat{A}}{n}} \| \Mz{Y}{n} - \Mn{\hat{A}}{n} \M{R}^\top \M{Q}_0^\top \|_F.
\end{equation}

After computing the Multi-TTM, we form $\M{W}_n = \Mz{Y}{n} \M{Q}_0$ 
to obtain the smaller least-squares problem
\begin{equation}\label{eq:ls-qr-q0-step}
\min_{\Mn{\hat{A}}{n}} \| \M{W}_n - \Mn{\hat{A}}{n} \M{R}^\top \|_F,
\end{equation}
and, finally, we use substitution with $\M{R}^\top$ to compute the factor matrix $\Mn{\hat{A}}{n}$. We call this algorithm CP-ALS-QR.

Another more stable way of solving \cref{eq:ls-qr-q0-step}, particularly when Khatri-Rao product $\M{Z}_n$ is rank deficient, is to use the singular value decomposition (SVD) of $\M{R}$ in addition to the QR factorization. Let $\M{R} = \M{U\Sigma V}^\top$ be the SVD. We then obtain 
\begin{equation*}
\min_{\Mn{\hat{A}}{n}}  \| \M{W} - \Mn{\hat{A}}{n} \M{V} \M{\Sigma} \M{U}^\top \|_F.
\end{equation*}
Now, since $\M{U}$ and $\M{V}$ are orthogonal and $\M{\Sigma}$ is a diagonal matrix, our solution to the least-squares problem is $\Mn{\hat{A}}{n} = \M{W}_n \M{U} \M{\Sigma}^{\dagger} \M{V}^\top$. Note that because we are utilizing the pseudoinverse on $\M{\Sigma}$, we can truncate small singular values and thereby manage a rank-deficient least-squares problem more stably. Using the SVD in this way gives us our second algorithm CP-ALS-QR-SVD. 

Both CP-ALS-QR and CP-ALS-QR-SVD are summarized in \cref{alg:cp-als-qr}, with a choice in \cref{line:qrsolve} to distinguish between the two methods.  
The two algorithms are implemented in MATLAB using the Tensor Toolbox \cite{tensor-toolbox}. 
The implementation is available here: \href{https://github.com/rlminste/CP-ALS-QR}{\textbf{https://github.com/rlminste/CP-ALS-QR}}. 

\begin{algorithm}
  \caption{CP-ALS-QR}
  \label{alg:cp-als-qr}
  \begin{algorithmic}[1]\footnotesize
    \Function{$[\bm{\lambda},\set{\Mn{A}{n}}]=$ CP-ALS-QR}{$\T{X},R$}\Comment{$\T{X}\in\mathbb{R}^{I_1\times \cdots \times I_N}$}
    \State Initialize factor matrices $\Mn{\hat{A}}{2}, \dots, \Mn{\hat{A}}{N}$ 
    \State \label{line:qr_factor}Compute compact QR-decomposition $\Mn{Q}{2} \Mn{R}{2}, \ldots, \Mn{Q}{N} \Mn{R}{N}$ of factor matrices
    \Repeat
    \For{$n=1,\dots, N$}
      \State \label{line:qr-khatrirao}$\M{V}_n \gets \Mn{R}{N} \Khat \cdots \Khat \Mn{R}{n+1} \Khat \Mn{R}{n-1} \Khat \cdots \Khat \Mn{R}{1}$ 
      \State \label{line:q0}Compute compact QR-decomposition $\M{V}_n = \M{Q}_0 \M{R}$ \Comment{Last step of QR decomposition}
      \State \label{line:ttm}$\T{Y} \gets \T{X} \times_1 \M{Q}_1^\top \times_2 \cdots \times_{n-1} \M{Q}_{n-1}^\top \times_{n+1} \M{Q}_{n+1}^\top \times_{n+2} \cdots \times_N \M{Q}_{N}^\top$\Comment{Multi-TTM}
      \State \label{line:apply_q0}$\M{W}_n \gets \Mz{Y}{n} \M{Q}_0$
      \State \label{line:qrsolve}Solve $\Mn{\hat{A}}{n} \M{R}^\top = \M{W}_n$ for $\Mn{\hat{A}}{n}$ by substitution (CP-ALS-QR) or SVD (CP-ALS-QR-SVD)
      \State Normalize columns of $\Mn{\hat{A}}{n}$ to obtain $\M{A}_n$ and $\V{\lambda}$
      \State \label{line:factorqr}Recompute QR-decomposition for updated factor matrix $\Mn{A}{n}=\Mn{Q}{n} \Mn{R}{n}$ 
    \EndFor
    \Until termination criteria met
    \State \textbf{return} $\bm{\lambda}$, factor matrices $\set{\Mn{A}{n}}$
    \EndFunction
  \end{algorithmic}
\end{algorithm}

\subsection{CP-ALS-QR Cost Analysis}\label{ssec:costs}

We now analyze the computational complexity of each iteration of CP-ALS-QR and CP-ALS-QR-SVD as presented in \cref{alg:cp-als-qr}.
Recall that $\T{X}$ has dimensions $I_1\times \cdots \times I_N$ and the CP approximation has rank $R$.
To simplify notation, we assume in the analysis that $I_1 \geq \cdots \geq I_N$.

The cost of forming the Khatri-Rao product $\M{V}_n$ of $N{-}1$ upper-triangular factors (\cref{line:qr-khatrirao}) is $R^N+\mathcal{O}(R^{N-1})$, assuming the individual Khatri-Rao products are formed pairwise and no sparsity is exploited.  
Here $\M{V}_n$ has dimensions $R^{N-1} \times R$, and it has $R^N/N+\mathcal{O}(R^{N-1})$ nonzeros.
We treat $\M{V}_n$ as a dense matrix here but discuss the possibility of exploiting sparsity in \cref{ssec:sparsity}.
Computing the QR decomposition of $\M{V}_n$ to obtain $\M{Q}_0$ and $\M{R}$ in \cref{line:q0} costs 
\begin{equation}
\label{eq:QRcost}
4R^{N+1}+\mathcal{O}(R^3),
\end{equation} 
assuming $\M{Q}_0$ is formed explicitly and again no sparsity is exploited.

The Multi-TTM is performed in \cref{line:ttm} and involves the input tensor.
We compute the resulting tensor $\T{Y}$, which has dimensions $R\times {\cdots} \times I_n\times {\cdots} \times R$, by performing single TTMs in sequence; to minimize flops we perform the $N{-}1$ TTMs in order of decreasing tensor dimension, which is left to right given our assumption above.
The cost of the first TTM is $2I_1\cdots I_NR$, the cost of the second is $2I_2\cdots I_N R^2$, and so on.
Thus, we can write the overall cost of the Multi-TTM as
\begin{equation}
\label{eq:TTMcost}
2I_1\cdots I_N R \left(1 + \frac{R}{I_1} + \frac{R^2}{I_1I_2} + \cdots + \frac{R^{N-2}}{I_1\cdots I_{n-1}I_{n+1}\cdots I_{N-1}} \right).
\end{equation} 

We apply $\M{Q}_0$ to $\Mz{Y}{n}$ via matrix multiplication in \cref{line:apply_q0} with cost $2I_nR^N$, assuming we use an explicit, dense $\M{Q}_0$.
Solving the linear system in \cref{line:qrsolve} costs $\mathcal{O}(I_nR^2)$, with an extra $\mathcal{O}(R^3)$ cost if the SVD of $\M{R}$ is computed for the CP-ALS-QR-SVD method.
(Note that using the more stable SVD approach to solve the linear system has no significant impact on the overall computational complexity.)
Finally, computing the QR decomposition of the updated $n$th factor matrix in \cref{line:factorqr} and forming the orthonormal factor $\M{Q}_n$ explicitly costs $4I_nR^2 + \mathcal{O}(R^3)$.

If the rank $R$ is significantly smaller than the tensor dimensions, then the cost is dominated by the first TTM, which has cost $2I_1\cdots I_NR$ from \cref{eq:TTMcost}.
In this case, the remaining TTMs are each at least a factor of $R/I_1$ times cheaper, and the computations involving $\M{Q}_0$ are  cheaper than any of the TTMs, because those computational costs are independent of any tensor dimensions.
The dominant cost of CP-ALS is the MTTKRP in \cref{line:mttkrp} of \cref{alg:cp-als}, which also has cost $2I_1\cdots I_NR$.
Thus, in the case of small $R$, the two algorithms have identical leading-order computational complexity per iteration.

If the rank $R$ is larger than all tensor dimensions, then the cost of the QR of $\M{V}_n$ given in \cref{eq:QRcost} will be the dominant cost.
If the rank is comparable to the tensor dimensions, then the computation and application of $\M{Q}_0$ in \cref{line:apply_q0}  and the subsequent TTMs after the first may also contribute to the running time in a significant way.

\subsection{Implementation Details and Extensions}\label{sec:implement}

\subsubsection{Efficient Computation of Approximation Error}\label{ssec:error}

For each of the CP algorithms (\cref{alg:cp-als,alg:cp-als-qr}), we consider computing the approximation error in two ways.
The more accurate but less efficient approach is to form the explicit representation of $\T{\hat{X}} = [\![ \V{\lambda}; \Mn{A}{1}, \ldots, \Mn{A}{N} ]\!]$ and compute the residual norm $\|\T{X} - \T{\hat{X}}\|$ directly.
The less accurate but more efficient approach exploits the identity $\|\T{X} - \T{\hat{X}}\|^2=\|\T{X}\|^2-2 \langle \T{X},\T{\hat{X}} \rangle + \|\T{\hat{X}}\|^2$ and computes $\langle \T{X},\T{\hat{X}} \rangle$ and $\|\T{\hat{X}}\|$ cheaply by using temporary quantities already computed by the ALS iterations ($\|\T{X}\|$ is pre-computed and does not change over iterations).

In the case of CP-ALS (\cref{alg:cp-als}), we have $\Mz{\hat{X}}{N} = \M{\hat A}_{N}\M{Z}_{N}^\top$, for $\Mn{\hat{A}}{N} = \Mn{A}{N} \cdot \text{diag}(\V{\lambda})$ and $\M{Z}_{N} = \M{A}_{N-1} \odot \dots \odot \M{A}_{1}$.
Then 
$$\langle \T{X}, \T{\hat{X}} \rangle = \langle \Mz{X}{N}, \M{\hat A}_{N}\M{Z}_{N}^\top \rangle = \langle \Mz{X}{N}\M{Z}_{N}, \M{\hat A}_N \rangle = \langle \M{M}_N,  \M{\hat A}_{N} \rangle,$$
where $\M{M}_N$ is the result of the MTTKRP computation in mode $N$, the mode of the last subiteration.
Thus, computing the inner product between the data and model tensors requires only $\mathcal{O} (I_NR)$ extra operations.
Likewise, we have
$$\|\T{\hat{X}}\|^2 = \langle \M{\hat A}_{N}\M{Z}_{N}^\top, \M{\hat A}_N \M{Z}_{N}^\top \rangle = \langle \M{Z}_{N}^\top \M{Z}_{N}, \M{\hat A}_{N}^\top \M{\hat A}_{N}\rangle = \langle \M{S}_N, \text{diag}(\V{\lambda})\M{G}_{N}\text{diag}(\V{\lambda}) \rangle,$$ 
where $\M{G}_{N} = \M{A}_{N}^\top \M{A}_{N}$ is the Gram matrix of the (normalized) $N$th factor and $\M{S}_N$ is the Hadamard product of the Gram matrices of the first $N-1$ modes.
Computing the norm of $\T{\hat{X}}$ thus requires only $\mathcal{O} (R^2)$ extra operations.
This efficient error computation is well known \cite{tensor-toolbox,EH+21,LK+17,TensorBox,SK16}.
Note that this approach is slightly less accurate than a direct computation: the identity applies to the square of the residual norm, so taking the square root of the difference of these quantities limits the accuracy of the relative error to the square root of machine precision.

We complete the efficient error computation for CP-ALS-QR (\cref{alg:cp-als-qr}) with similar cost.
In this case, we have $\Mz{\hat{X}}{N} = \M{\hat A}_{N} \M{Z}_{N}^\top$ with $\M{Z}_{N} = (\M{Q}_{N-1} \Kron \cdots \Kron \M{Q}_{1}) \M{Q}_0 \M{R}$.
Then 
$$\langle \T{X}, \T{\hat{X}} \rangle = \langle \Mz{X}{N}, \M{\hat A}_{N}\M{Z}_{N}^\top \rangle = \langle \Mz{X}{N}(\M{Q}_{N-1} \Kron \cdots \Kron \M{Q}_{1}) \M{Q}_0, \M{\hat A}_{N}\M{R}^{\top} \rangle = \langle \M{W}_N, \M{\hat A}_{N}\M{R}^{\top} \rangle,$$
where in \cref{alg:cp-als-qr}, $\M{W}_N$ is the result of the Multi-TTM (except in mode $N$) and the multiplication with $\M{Q}_0$, and thus has dimension $I_N \times R$.
Thus, the cost of this computation is dominated by that of computing $\M{\hat A}^{(N)}\M{R}^{\top}$, or $\mathcal{O} (I_NR^2)$.
We also have
$$\|\T{\hat{X}}\|^2 = \langle \M{\hat A}_{N}\M{Z}_{N}^\top, \M{\hat A}_{N}\M{Z}_{N}^\top \rangle = \langle \M{Z}_{N}^\top \M{Z}_{N}, \M{\hat A}_{N}^\top \M{\hat A}_{N}\rangle = \langle \M{R}^\top \M{R}, \text{diag}(\V{\lambda})\M{R}_{N}^\top \M{R}_{N}\text{diag}(\V{\lambda}) \rangle,$$
where $\M{R}_{N}$ is the triangular factor in the QR decomposition of $\M{A}_N$.
The cost of this extra computation is $\mathcal{O} (R^3)$.

\subsubsection{Kruskal Tensor Input}\label{ssec:ktensor}
The analysis of both CP-ALS-QR and CP-ALS-QR-SVD, as explained in \cref{ssec:costs}, assume the input tensor is a dense tensor. 
When the input tensor has special structure, the key operations can be computed more efficiently.
We also implemented a version of each algorithm that exploits inputs with Kruskal structure, that is, a tensor stored as factor matrices and corresponding weights, which we use for the input in \cref{ssec:sine}. 
Exploiting this structure is beneficial as we avoid forming the input tensor or Multi-TTM product tensor $\T{Y}$, because all computations can be performed using the factor matrices instead. 

Note that the Tensor Toolbox has optimized the MTTRKP (\cref{alg:cp-als}, \cref{line:mttkrp}) and Multi-TTM (\cref{alg:cp-als-qr}, \cref{line:ttm}) computations for a Kruskal tensor input \cite{efficient-matlab}. For example, the complexity of the Multi-TTM product to obtain \cref{eq:ls-qr-ttm-step} in the case of a Kruskal tensor becomes $\mathcal{O} (R^N \sum_{i=2}^N I_i)$. For an $N$-mode Kruskal tensor $\T{X} \in \mathbb{R}^{I_1 \times \cdots \times I_N}$ of rank $R$ and $N$ matrices $\M{V}_j \in \mb{R}^{I_j \times S}$ for $j = 1,\dots, N$, the MTTKRP
\begin{equation*}
\Mz{X}{n} ( \Mn{V}{N} \Khat \cdots \Khat \Mn{V}{n+1} \Khat \Mn{V}{n-1} \Khat \cdots \Khat \Mn{V}{1})
\end{equation*}
has cost $\mathcal{O} (RS \sum_{i=1}^N I_i)$. This is an improvement compared to the cost of the MTTKRP in the dense case, which is on the order of the product of the $N$ dimensions instead of their sum.

The Kruskal structure has the added benefit of reducing the cost of computing the product $\M{W}_n = \Mz{Y}{n}\M{Q}_0$ in our QR-based algorithms (see \cref{line:apply_q0} in \cref{alg:cp-als-qr}), as we do not need to explicitly matricize $\T{Y}$. 
Let us consider the first mode as an example. 
In this case, $\T{Y}$ is an $N$-dimensional Kruskal tensor with rank $R$ and $\T{Y} =  [\![ \lambda; \Mn{B}{1}, \Mn{B}{2}, \ldots ,\Mn{B}{N} ]\!]$, with $\Mn{B}{1} \in \mb{R}^{I_1 \times R}$, $\Mn{B}{i} \in \mb{R}^{R \times R}$  for $i = 2, \ldots, N$, and $\M{Q}_0 \in \mb{R}^{R^{N-1} \times R}$. 
Then, we can write $\Mz{Y}{1} = \Mn{B}{1} ( \Mn{B}{N} \Khat \cdots \Khat \Mn{B}{2})^\top$, treat $\M{Q}_0$ as a matricized tensor and compute the $\M{W} = \Mn{B}{1} [( \Mn{B}{N} \Khat \cdots \Khat \Mn{B}{2})^\top \M{Q}_0 ]$ using an MTTKRP followed by a small matrix product. 
This gives us a total cost of $2(I_1 R^2 + R^{N+1})$. 

\subsubsection{Other Computation-Reducing Optimizations}
\label{ssec:sparsity}

As noted in \cref{ssec:costs}, the Khatri-Rao product of triangular matrices $\M{V}_n$ computed in \cref{alg:cp-als-qr} is a sparse matrix with density proportional to $1/N$, where $N$ is the number of modes.
This is because the $i$th column of $\M{V}_n$ is a Kronecker product of $N{-}1$ vectors each with $i$ nonzeros and therefore has $i^{N-1}$ nonzeros.
This sparsity could be exploited in the computation of $\M{V}_n$ (\cref{line:qr-khatrirao}), computing its QR decomposition (\cref{line:q0}), and applying its orthonormal factor $\M{Q}_0$ (\cref{line:apply_q0}).
In particular, when using a sparse QR decomposition algorithm, there will be no fill-in, as every row is dense to the right of its first nonzero, and the number of flops required is a factor of $\mathcal{O}(1/N^2)$ times that of the dense QR algorithm.
The computational savings in computing $\M{V}_n$ and applying $\M{Q}_0$ is $\mathcal{O}(1/N)$.
However, the use of (general) sparse computational kernels comes at a price of performance, so for small $N$ we do not expect much reduction in time and did not exploit sparsity in our implementation.

An important optimization for CP-ALS (\cref{alg:cp-als}) is to avoid recomputation across MTTKRPs of the different modes.
For example, the Khatri-Rao products $\M{Z}_n$ and $\M{Z}_{n+1}$ share $N-2$ different factors, so the computations of $\M{M}_n$ and $\M{M}_{n+1}$ have significant overlap. 
The general approach to avoiding this recomputation is known as dimension trees, as a tree of temporary matrices can be computed, stored, and re-used for the MTTKRPs across modes \cite{EH+21,PTC13a}.
Using dimension trees reduces the outer-iteration CP-ALS cost from $N$ MTTKRPs to the cost of 2 MTTKRPs.

Similar savings can be obtained by applying dimension trees to the set of Multi-TTM operations in CP-ALS-QR (\cref{alg:cp-als-qr}).
In this case, we exploit the overlap in individual TTMs across modes and store a different set of intermediate tensors that can be re-used across modes.
Dimension trees have been used for Multi-TTM before in the context of Tucker decompositions for sparse tensors and the Higher-Order Orthogonal Iteration algorithm \cite{BMVL12,KU16}.
As in the case of CP-ALS, dimension trees can reduce the outer-iteration CP-ALS-QR cost from $N$ TTMs involving the data tensor to 2 TTMs.
Neither of these reductions come at the expense of lower performance, so we can expect $\mathcal{O}(N)$ speedup in each case.
For CP-ALS-QR, there are other overlapping computations that can be similarly exploited.
For example, the QR decomposition of $\M{V}_n$ can be performed using a tree across the Khatri-Rao factors, some of which are shared across modes, though the structure of the orthonormal factor would need to be maintained when applying it.
Because the Tensor Toolbox implementation of CP-ALS does not employ dimension trees, for fair comparison, we do not use them for CP-ALS-QR either.


\section{Numerical Experiments}\label{sec:examples}
In this section, we explore several examples that demonstrate the benefits of CP-ALS-QR and CP-ALS-QR-SVD over the typical ALS approaches. Specifically, we will demonstrate the performance of our algorithms as well as show the stability of our algorithms by considering ill-conditioned problems.

\subsection{Performance Results}
\label{ssec:perf}
As seen in our analysis in \cref{ssec:costs}, the dominant cost for our new algorithms is the same as for CP-ALS and CP-ALS-PINV when $R$ is small. 
For large $R$, we see that the lower-order terms for CP-ALS-QR and CP-ALS-QR-SVD do have an effect on the runtime. 
We verify the comparable runtimes for small $R$ and examine the slowdown for large $R$ with a few experiments here.

We break each algorithm down to its key components and time each individually. 
These components are listed in \cref{tab:its_parts}. 
Each row of the table represents corresponding parts of the two different types of algorithm.
\begin{table}[!ht]
\centering
\begin{tabular}{|c|c|}
	\hline
	\textbf{CP-ALS, CP-ALS-PINV} & \textbf{CP-ALS-QR, CP-ALS-QR-SVD} \\
	\hline
	MTTKRP & Multi-TTM \\
	Gram of factor matrices & QR of factor matrices \\
	 N/A & Computing $\M{Q}_0$ \\
	 N/A & Applying $\M{Q}_0$\\
	 Other & Other \\
	 \hline
\end{tabular}
\caption{Breakdown of main components in each iteration of CP-ALS, CP-ALS-PINV, CP-ALS-QR, and CP-ALS-QR-SVD. We use this breakdown in our performance experiments shown in \cref{fig:its}.}
\label{tab:its_parts}
\end{table}
For CP-ALS and CP-ALS-PINV (\cref{alg:cp-als}), the MTTKRP refers to forming $\M{M}_n = \Mz{X}{n} \M{Z}_n$, see \cref{line:mttkrp}, while we compute the Gram matrices for each factor matrix in \cref{line:gram}. For CP-ALS-QR and CP-ALS-QR-SVD (\cref{alg:cp-als-qr}), the Multi-TTM is computed when applying the Kronecker product of $\M{Q}_j$ matrices to $\Mz{X}{n}$, see \cref{line:ttm}, and we compute the QR factorization of each factor matrix in \cref{line:qr_factor}. Computing $\M{Q}_0$ involves computing a QR factorization of Khatri-Rao product $\M{V}_n$, see \cref{line:q0}, and applying $\M{Q}_0$ is a matrix multiplication, see \cref{line:apply_q0}. The steps included in the ``Other'' category include solving (by substitution or SVD), finding our weight vector $\V{\lambda}$, and the error computation. These steps are combined as none represent a significant portion of the runtime for any of the four algorithms. 

\begin{figure}[!t]
\centering
\includegraphics[scale=.4]{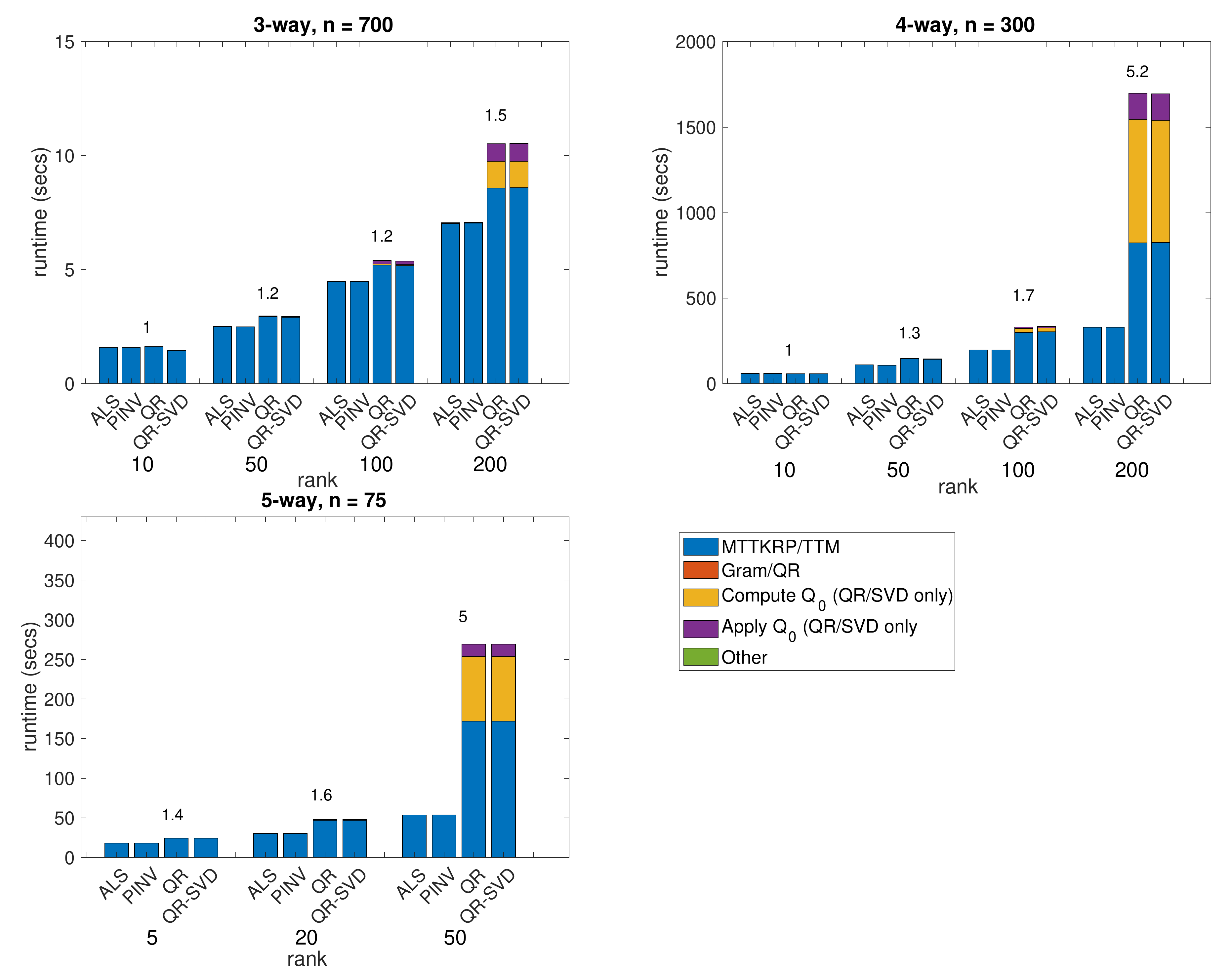}
\caption{Average runtime in seconds of a single iteration for CP-ALS, CP-ALS-PINV, CP-ALS-QR, and CP-ALS-QR-SVD for a three-way tensor of size 700 (top left), a four-way tensor of size 300 (top right), and a five-way tensor of size 75 (bottom left). Results are plotted for increasing rank values, and the slowdown ratio between the runtimes of CP-ALS-QR and CP-ALS is plotted above the group of results for each rank.}
\label{fig:its}
\end{figure}

The three tensors we test are randomly generated cubical tensors of three, four, and five modes. The three-way tensor has dimension $700$, the four-way tensor has dimension $300$, and the five-way tensor has dimension $75$. We computed the average iteration time over 10 iterations (omitting the first iteration to ensure a warm cache). The tolerance we used for all algorithms was $10^{-10}$, and we computed the error in the efficient manner described in \cref{ssec:error}.
The results for these three tensors with increasing rank values are in \cref{fig:its}. For each rank value, we also plot the slowdown ratio we see between the overall runtimes of CP-ALS-QR and CP-ALS.

For all three tensors, the dominant cost per iteration is the MTTKRP for CP-ALS and CP-ALS-PINV, and the Multi-TTM for CP-ALS-QR and CP-ALS-QR-SVD. Only for high ranks do other costs, computing and applying the QR of the Khatri-Rao product, even appear visibly in the plot. In the three-way case, all the slowdown ratios are close to $1\times$. This is similar for low ranks in four and five modes, but the ratio for high ranks jumps up to $5\times$.
These results demonstrate that the slowdown incurred by using our QR-based algorithms is not significant when the CP rank is small.

\subsection{Collinear Factor Matrices}
\label{ssec:collinear}
In this example, we test our algorithms on synthetic tensors constructed so that the factor matrices are ill-conditioned. Following the approach in \cite{score-info}, we create this tensor from randomly generated factor matrices and weights so that we are able to compare the results of our algorithms to the true solution, and we add Gaussian noise. 
The randomly generated factor matrices are constructed as in \cite{tomasi2006comparison} so that we can control their collinearity.
For our experiments, we construct a 3-way $50 \times 50 \times 50$ tensor with rank 5 and varying levels of collinearity and noise.
We compare our QR-based algorithms, CP-ALS-QR and CP-ALS-QR-SVD, to CP-ALS and CP-ALS-PINV. We also compare all four ALS algorithms to an optimization-based CP method using Gauss-Newton implemented in Tensorlab \cite{tensorlab} we describe in \cref{ssec:gn}. 
We test all combinations of three different noise levels $10^{-4}, 10^{-7}$, and $10^{-10}$, and three different collinearity levels $1-10^{-4}, 1-10^{-7}$, and $1-10^{-10}$. For each configuration, we run $100$ trials of each algorithm to approximate a rank-$5$ CP factorization. The maximum number of iterations is $500$ and the convergence tolerance for change in the relative error is $10^{-15}$. 
We use such a tight tolerance to ensure the computed metrics reflect what is attainable by the algorithm and not an artifact of early convergence.
A random guess is used for each initialization, and each algorithm is configured with the same initial factor matrices.

To measure the performance of each algorithm, we consider the number of iterations to converge, the relative error defined as $\||\T{X} - \T{\hat{X}}\|/\|\T{X}\|,$ and a similarity measure between the approximation and true tensor called the score \cite{tomasi2006comparison}. 
The score is a proxy for forward error by checking if two given Kruskal tensors are nearly equivalent up to scaling and permutation. 
The maximum score is 1 (for Kruskal tensors that are equivalent), and the minimum is 0.
An alternative metric called CorrIndex \cite{CorrIndex} also measures how closely two Kruskal tensors match. 
We observed similar quantitative behaviors between the score and CorrIndex and report only the score in our experiments.

Results of the experiments are presented in \cref{tab:coll-noise}, where we see that the combination of noise and collinearity affects the ill-conditioning of the problem in different ways. When the noise level is $10^{-4}$ (first row), as well as when the collinearity is $1-10^{-4}$ (first column), all ALS algorithms have similar performances, implying that the subproblems are not ill-conditioned in these cases. 
The behavior of Gauss-Newton is more variable across the first row and down the first column.
In the first row ($10^{-4}$ noise), we see the score decrease dramatically as collinearity increases.
In the first column ($1-10^{-4}$ collinearity), we see much faster (quadratic) convergence as noise decreases, and Gauss-Newton converges to the correct solution.  

We observe variability across ALS algorithms in the bottom-right $2\times 2$ grid of experiments, where the combination of higher collinearity ($1-10^{-7}$ and $1-10^{-10}$) and low noise ($10^{-7}$ and $10^{-10}$) creates ill-conditioned subproblems.
Our first observation is that CP-ALS-QR and CP-ALS-QR-SVD are robust in these cases: they obtain the lowest relative error with very little variation across initialization, they converge in fewer iterations than CP-ALS and CP-ALS-PINV, and their scores are generally high.
Among the two normal equations based algorithms, we observe that CP-ALS-PINV obtains a better score than CP-ALS in all cases, but this comes at the expense of higher backward error.
These two algorithms rarely converge before they hit the maximum of 500 iterations, which we attribute to the instability of the iterations, particularly near the ill-conditioned solution.
The behavior of Gauss-Newton is consistent across these four cases: it converges quickly (due to the low noise), but the scores and relative errors are generally worse than the ALS-based algorithms, sometimes significantly.
Comparing CP-ALS-QR and CP-ALS-QR-SVD, we see that the SVD variant is slightly more robust, obtaining higher scores, but sometimes requires more iterations to converge.


\begin{table}
    \centering
    \begin{tabular}{|c|M{49mm}M{46mm}M{45mm}|}
       \hline  
        {} & \multicolumn{3}{c|}{}\\[1pt]
        {} & \multicolumn{3}{c|}{\small Collinearity}\\[8pt]
       \hline 
        {} &  \multicolumn{3}{c|}{} \\[1pt]
       \small Noise & \footnotesize$1-10^{-4}$ & \footnotesize$1-10^{-7}$ & \footnotesize$1-10^{-10}$  \\[8pt]
       \hline
       \footnotesize$10^{-4}$ & \includegraphics[scale=.25]{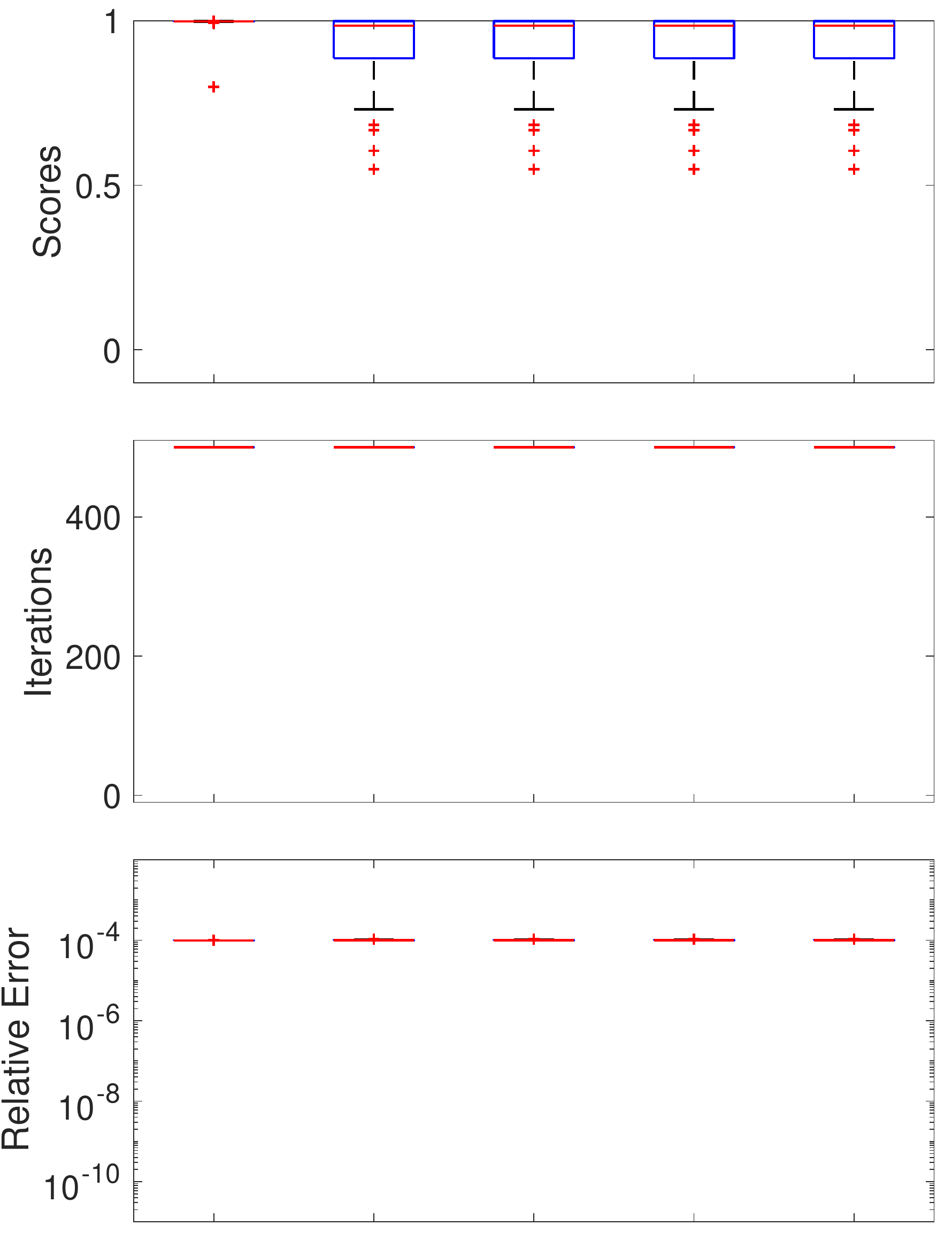} & \includegraphics[scale=.25]{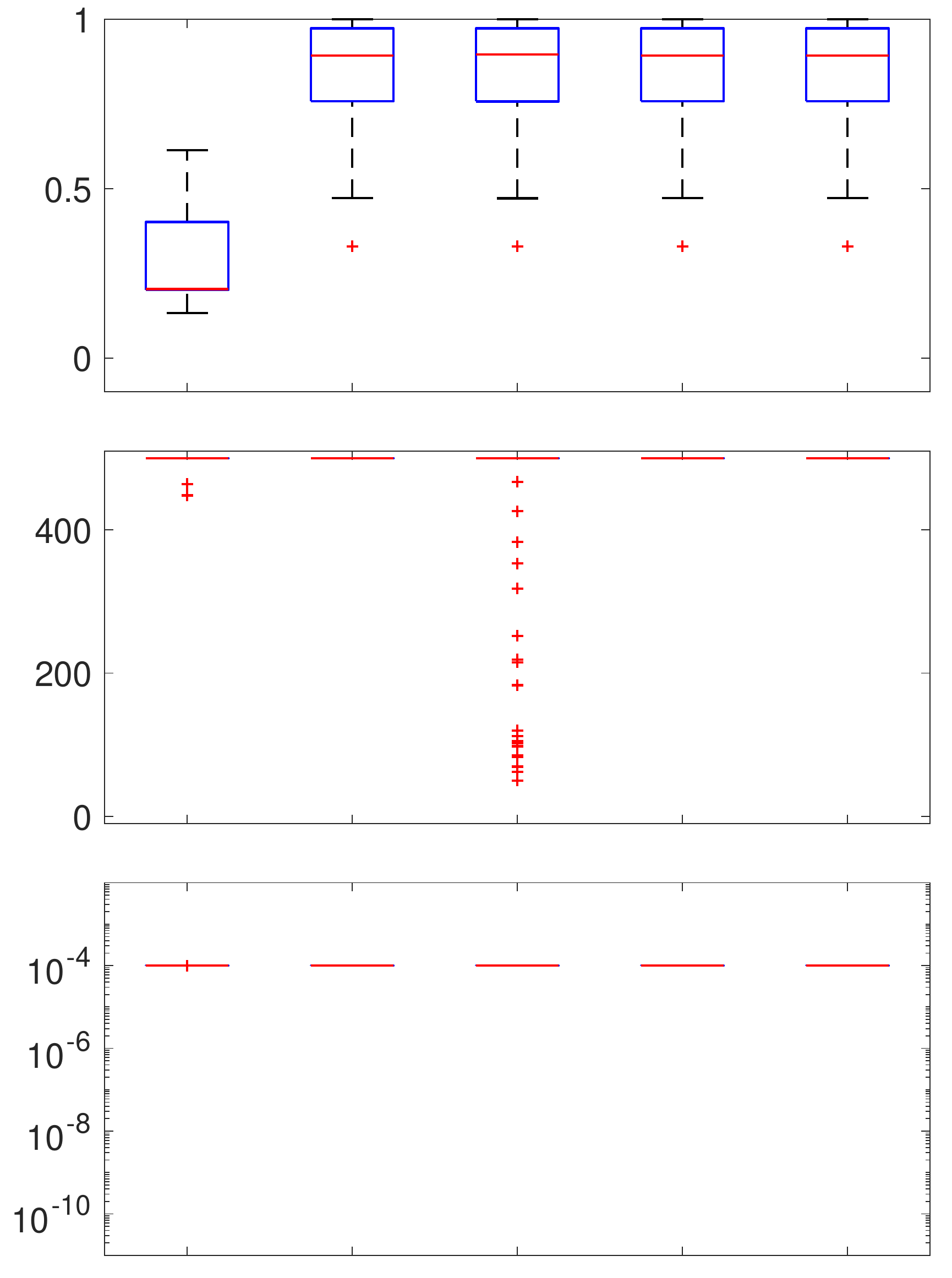} & \includegraphics[scale=.25]{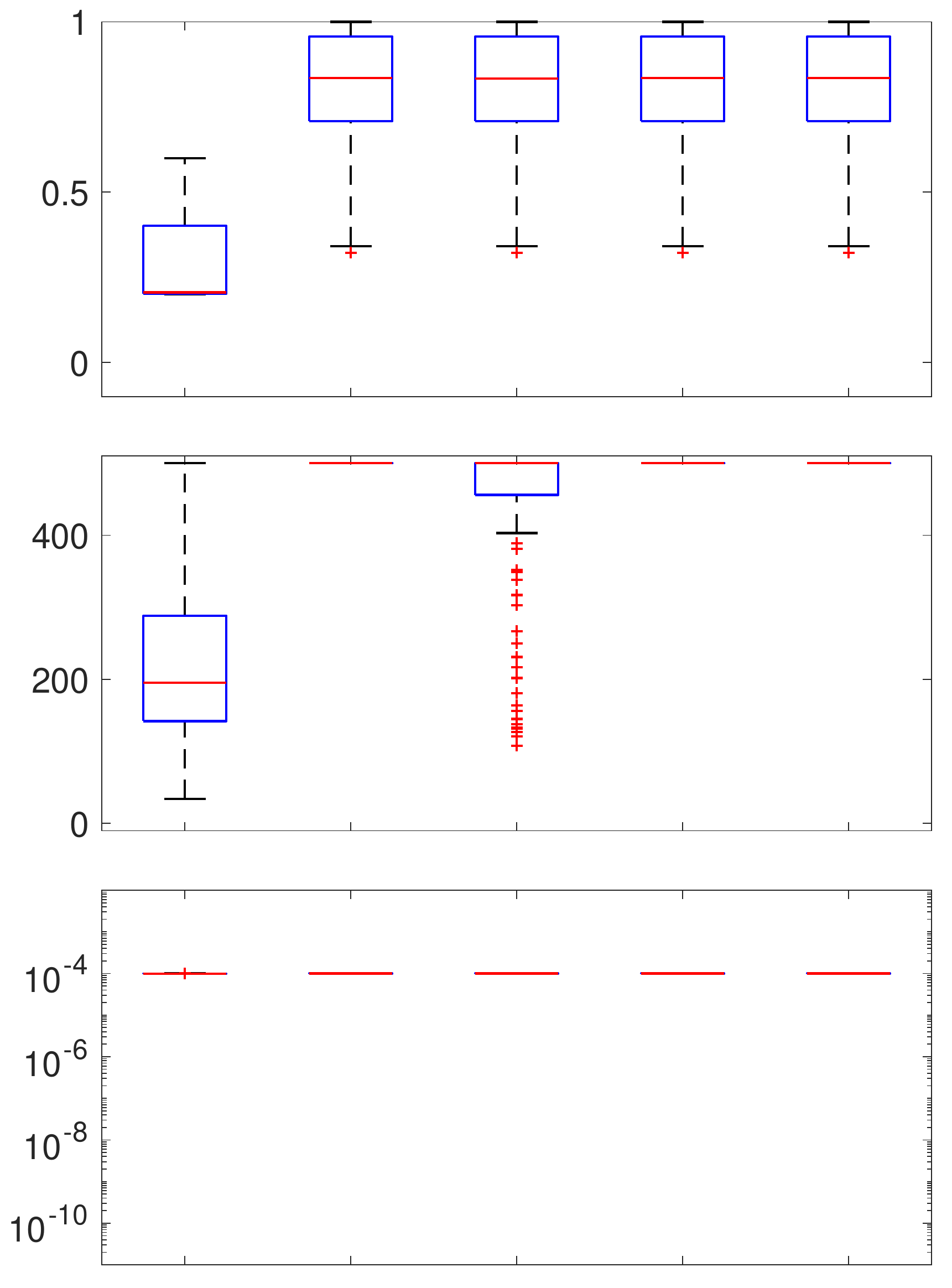}\\
	\hline
        \footnotesize$10^{-7}$ & \includegraphics[scale=.25]{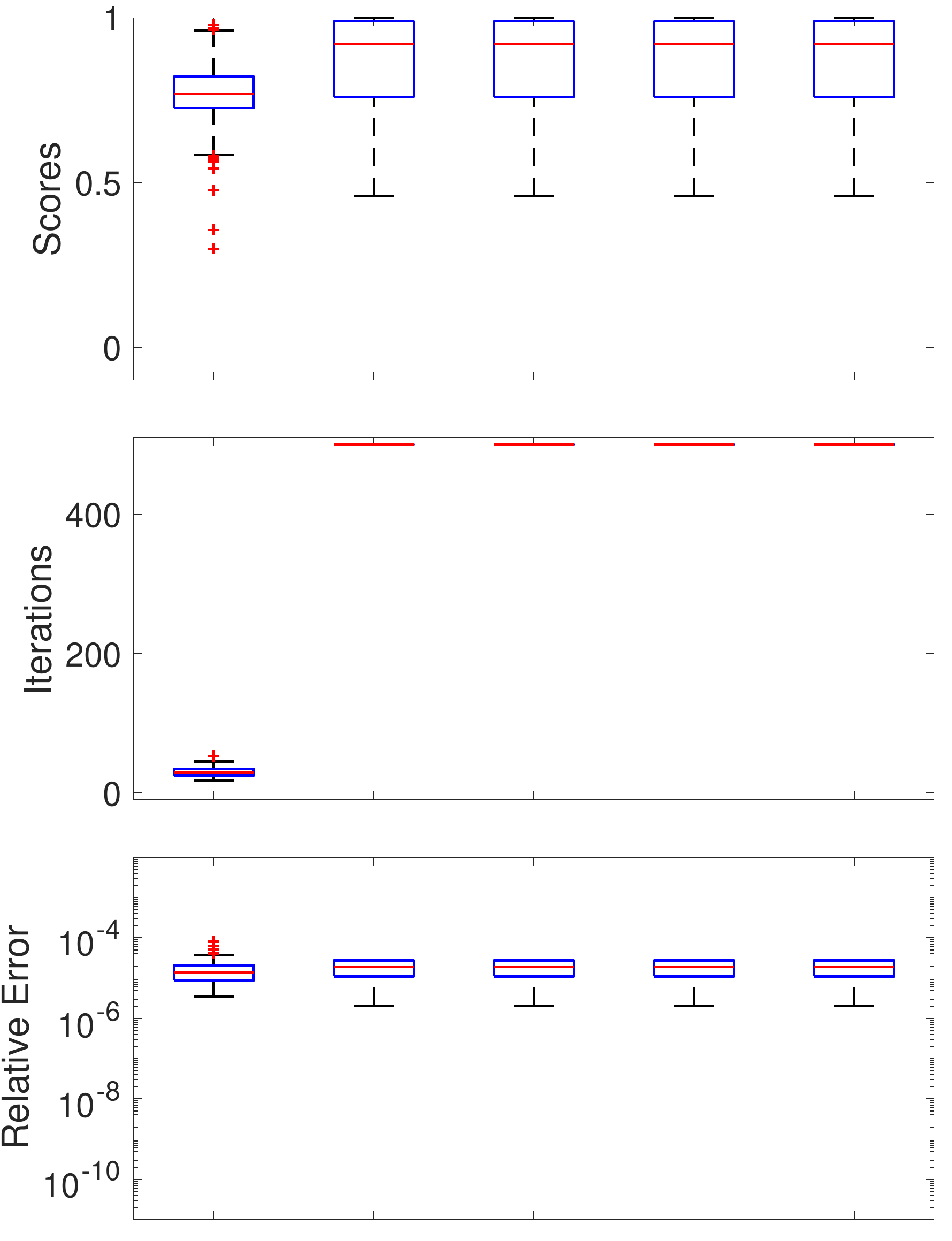} & \includegraphics[scale=.25]{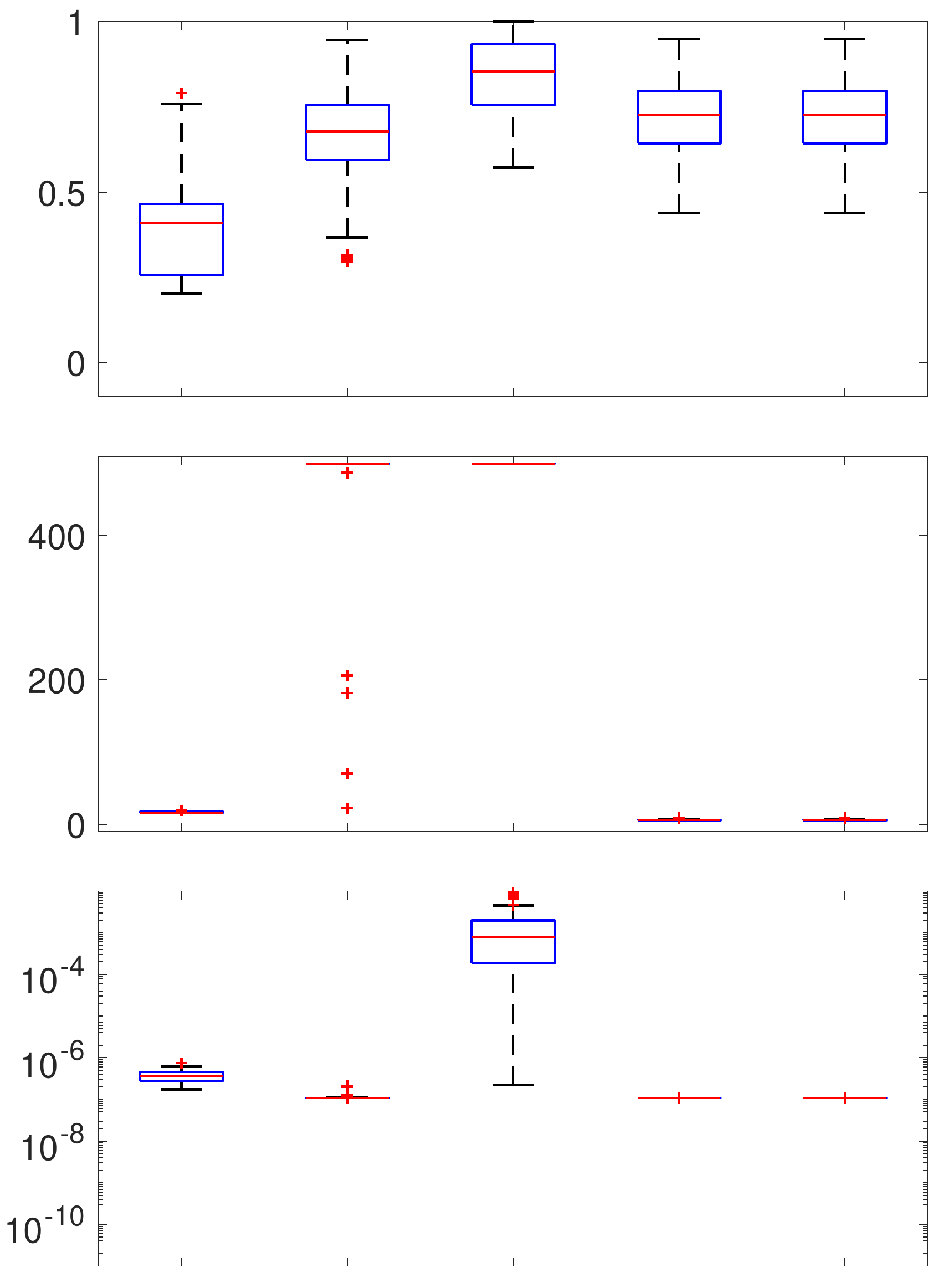} & \includegraphics[scale=.25]{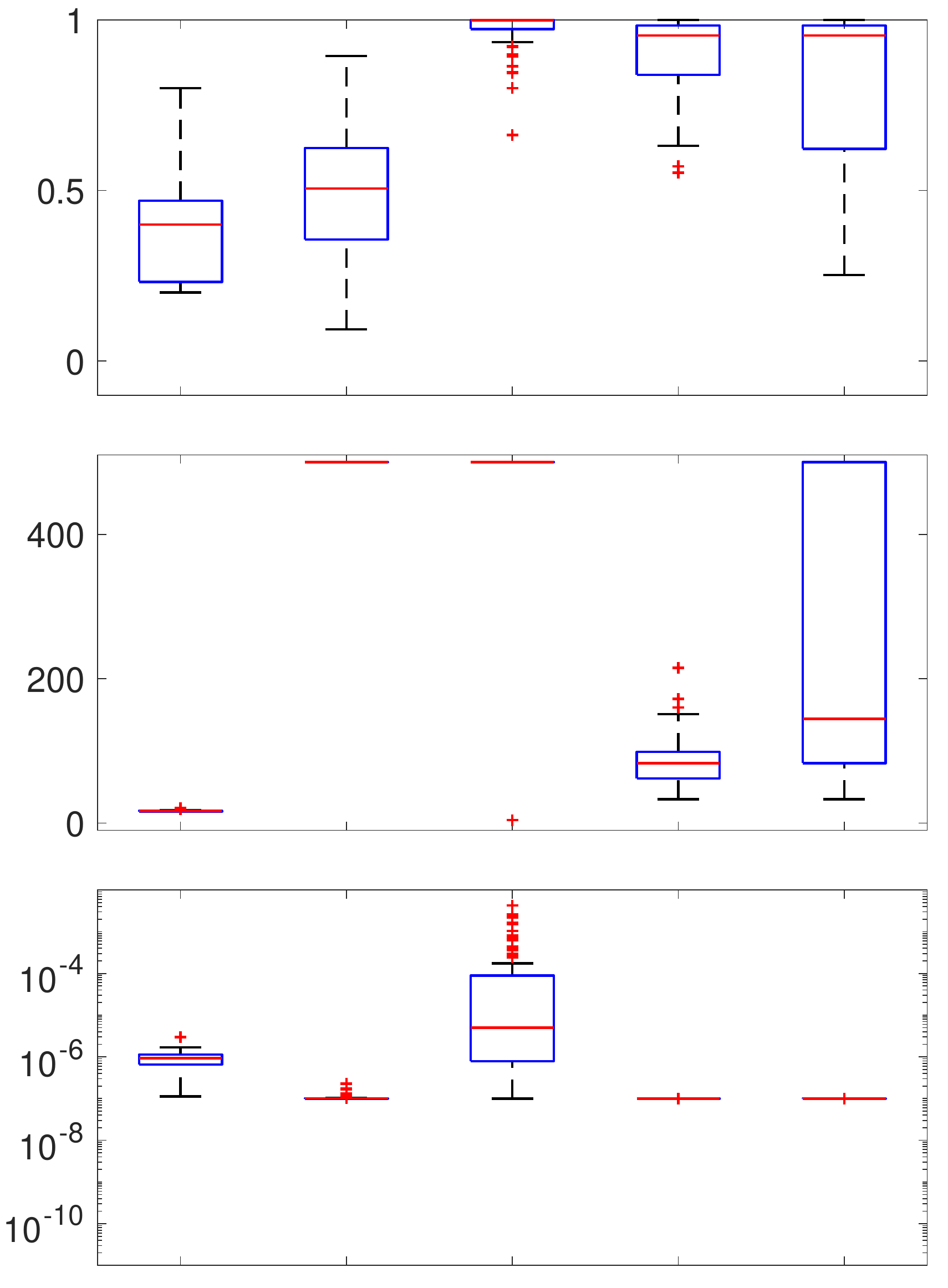}\\
        \hline
        	 \footnotesize$10^{-10}$ & \hspace{-.4cm} \includegraphics[scale=.25]{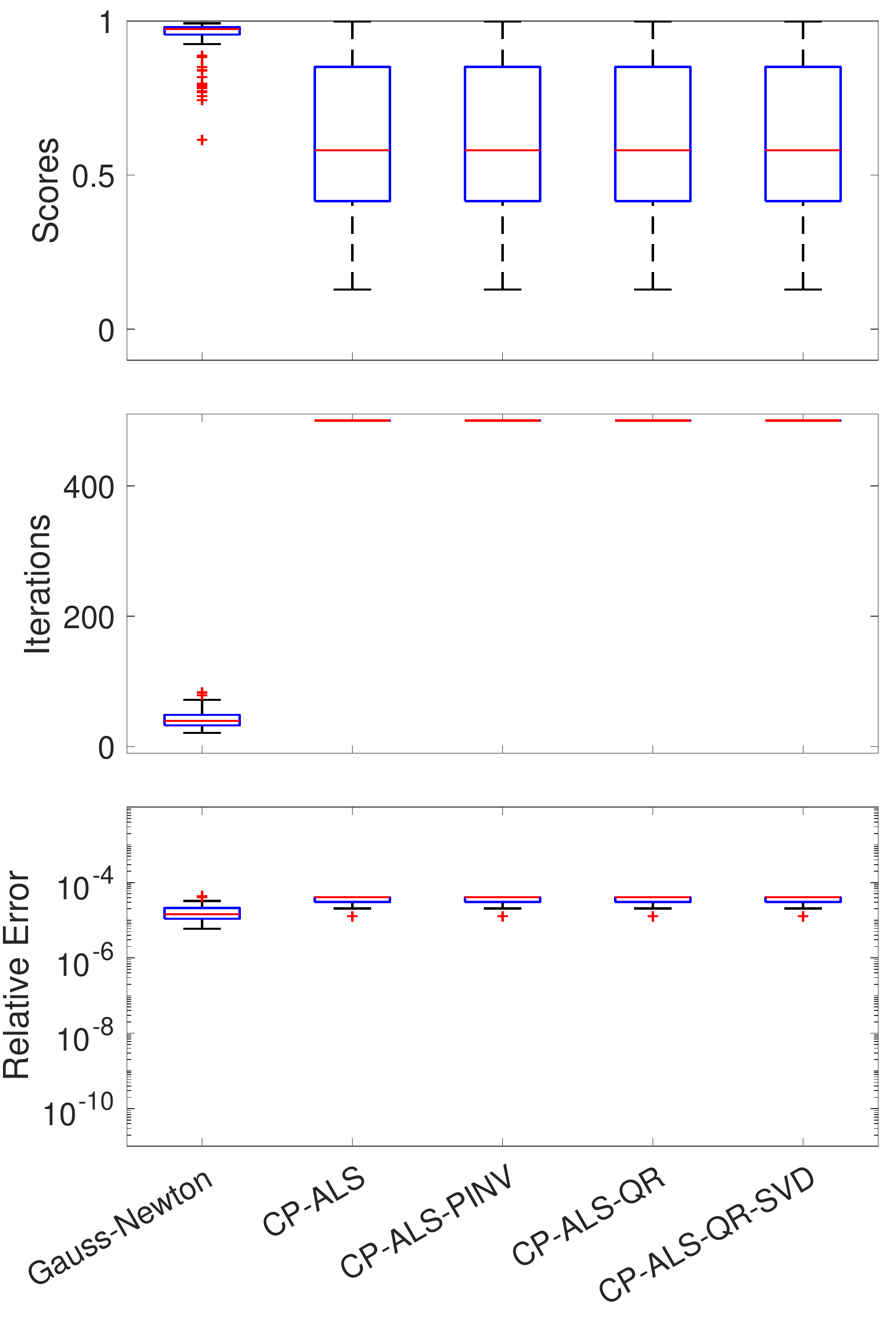}  & \includegraphics[scale=.25]{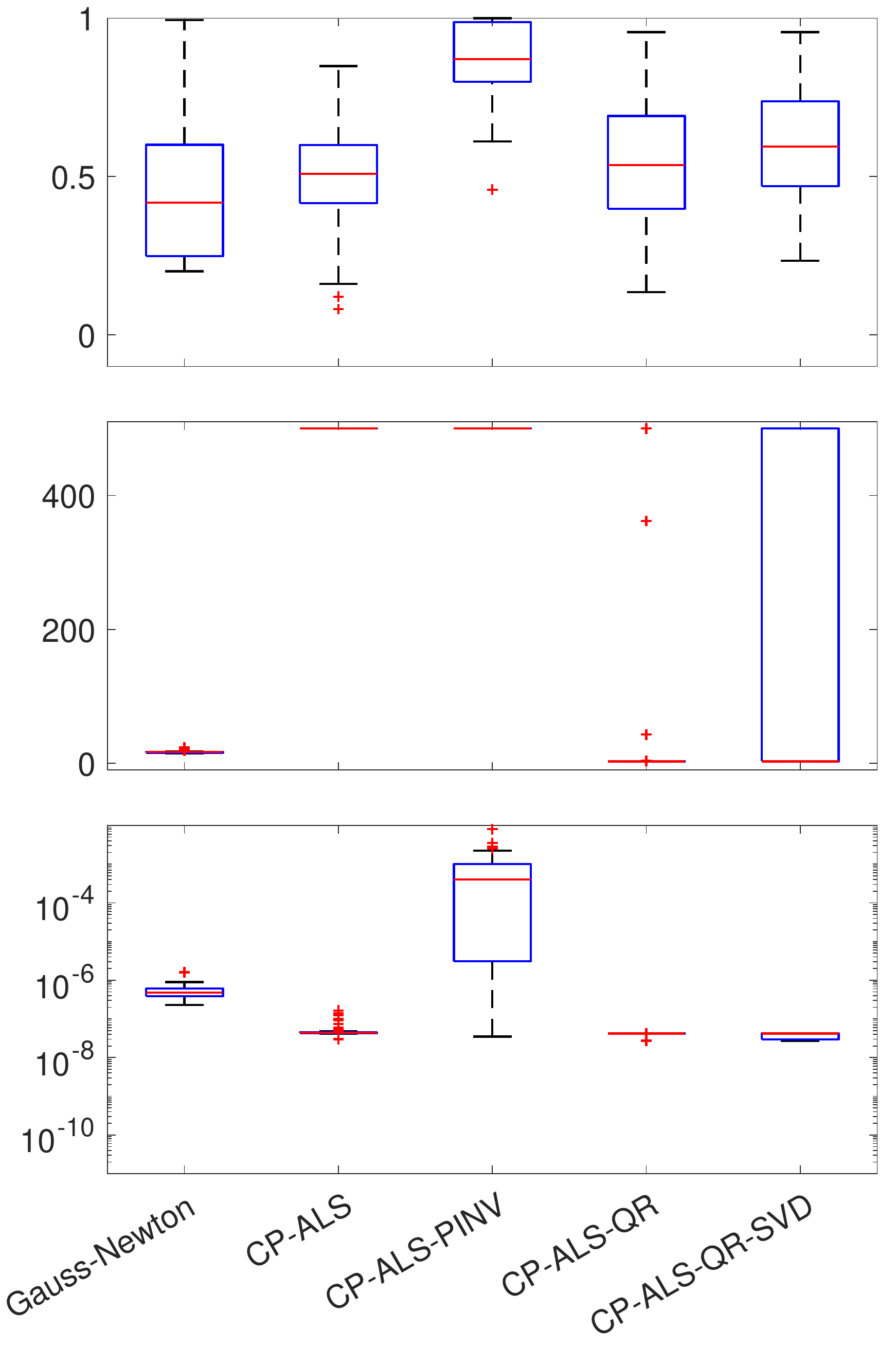} & \includegraphics[scale=.25]{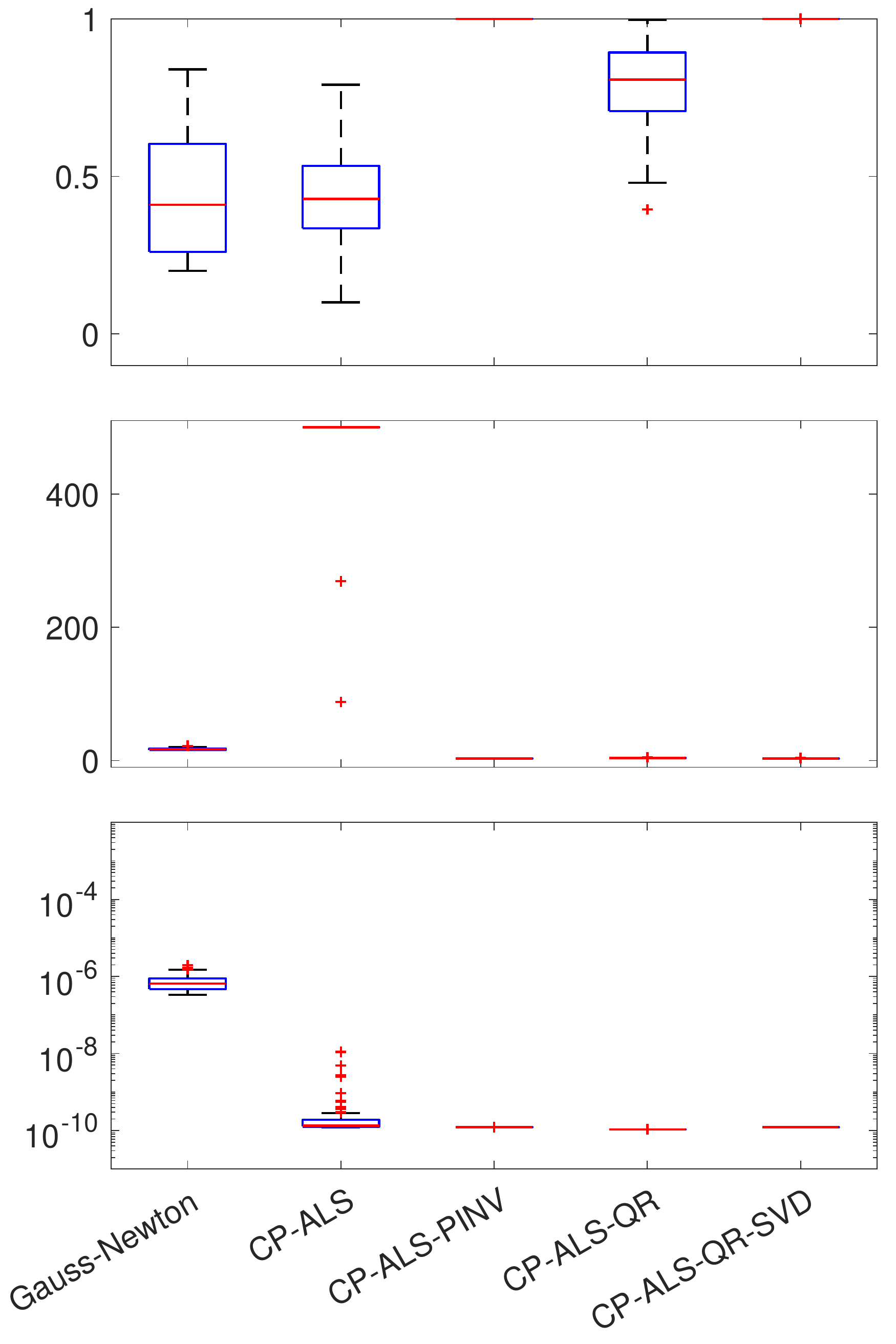} \\
	 \hline
    \end{tabular}
    \caption{Scores, iterations, and relative error boxplots for Gauss-Newton, CP-ALS, CP-ALS-PINV, CP-ALS-QR, and CP-ALS-QR-SVD on a $50 \times 50 \times 50$ synthetic tensor of rank 5 with three different levels of collinearity for the true factor matrices and three different levels of Gaussian noise added.}
    \label{tab:coll-noise}
\end{table}

To summarize, when subproblems are well conditioned (either high noise or low collinearity), we observe no difference in the convergence or accuracy in the ALS algorithms.
In the presence of ill-conditioned subproblems, we see that the QR-based algorithms, CP-ALS-QR and CP-ALS-QR-SVD, have stable performances in all scenarios, while the normal equations based algorithms tend not to converge quickly and also suffer from higher forward or backward errors.
We see that using the SVD within CP-ALS-PINV can improve the quality of solution compared to CP-ALS, it sacrifices the backward error and is not as robust as the QR-based methods.
Gauss-Newton is a viable alternative to ALS, exhibiting much faster convergence in the case of low noise, but it suffers from the same sensitivity to ill-conditioned problems and is not as robust as CP-ALS-QR or CP-ALS-QR-SVD in those cases.

\subsection{Sine of Sums}
\label{ssec:sine}
We further test the stability of our QR-based algorithms on a function approximation problem discussed in \cite{sine-defn}. The sine of sums $\sin(x_1 + \dots + x_N)$ is an example of an $N$-dimensional function that can be approximated in such a way where complexity grows linearly with $N$ instead of exponentially. 
These efficient approximations are sometimes referred to as separated representations, and they are closely related to CP decompositions. 
In this case, we are simulating the numerical discovery of an efficient separated representation, as we already know it exists.
That is, given a separated representation with large rank as input, we seek to compute a lower-rank representation that approximates it to numerical precision.
We consider this problem as it can be ill-conditioned depending on the representation, which is nonunique.

\subsubsection{Setup}
The multivariate sine of sums function $\sin(x_1 + \dots + x_N)$ can be discretized as a dense $N$-mode tensor $\T{X} \in \mb{R}^{n \times \dots \times n}$ with $x_j \in \mb{R}^{n}$ as vectors discretizing the interval $[0,2\pi)$ for $j = 1,\ldots,N$. 
As the sine of sums can be expressed as a sum of $2^{N-1}$ terms using sum identities for sine and cosine, we can expand all terms to obtain a rank-$2^{N-1}$ representation of $\T{X}$, which corresponds to an exact CP decomposition.  
Another exact CP representation of rank $N$ also exists, of the form 
\begin{equation}\label{eq:sinsum_rkd}
	\T{X} = \sin\left(\sum_{j=1}^N x_j \right) = \sum_{j=1}^N \sin(x_j) \prod_{k=1,k\neq j}^N \frac{\sin(x_k+\alpha_k - \alpha_j)}{\sin(\alpha_k - \alpha_j)},
\end{equation}
where $\alpha_j$ must satisfy $\sin(\alpha_k - \alpha_j) \neq 0$ for all $j \neq k$. This rank-$N$ representation is nonunique, and can be numerically unstable depending on the choice of $\alpha_j$. 

The input representations are already rank-$2^{N-1}$ Kruskal tensors, and we exploit that structure in our algorithms as discussed in \cref{ssec:ktensor}. 
In the experiments below, we vary the number of modes $N$, corresponding to the number of variables in the sine of sums function, and the dimension $n$ of each mode, corresponding to the number of discretization points in the interval $[0,2\pi)$.
In each experiment, we consider the relative error of four algorithms, CP-ALS, CP-ALS-PINV, CP-ALS-QR, and CP-ALS-QR-SVD, at the end of each of the first 40 iterations. 
We use the same random initial guesses for the factor matrices across algorithms and a convergence tolerance of 0. 
We also compute the relative error $\| \T{X} - \T{\hat{X}} \|/ \|\T{X}\|$ directly (as opposed to the methods described in \cref{ssec:error}) as we need a more accurate computation of the error to truly compare the accuracy of our algorithms.

\subsubsection{Lower-order tensors}
\begin{figure}
\centering
\includegraphics[scale=.4]{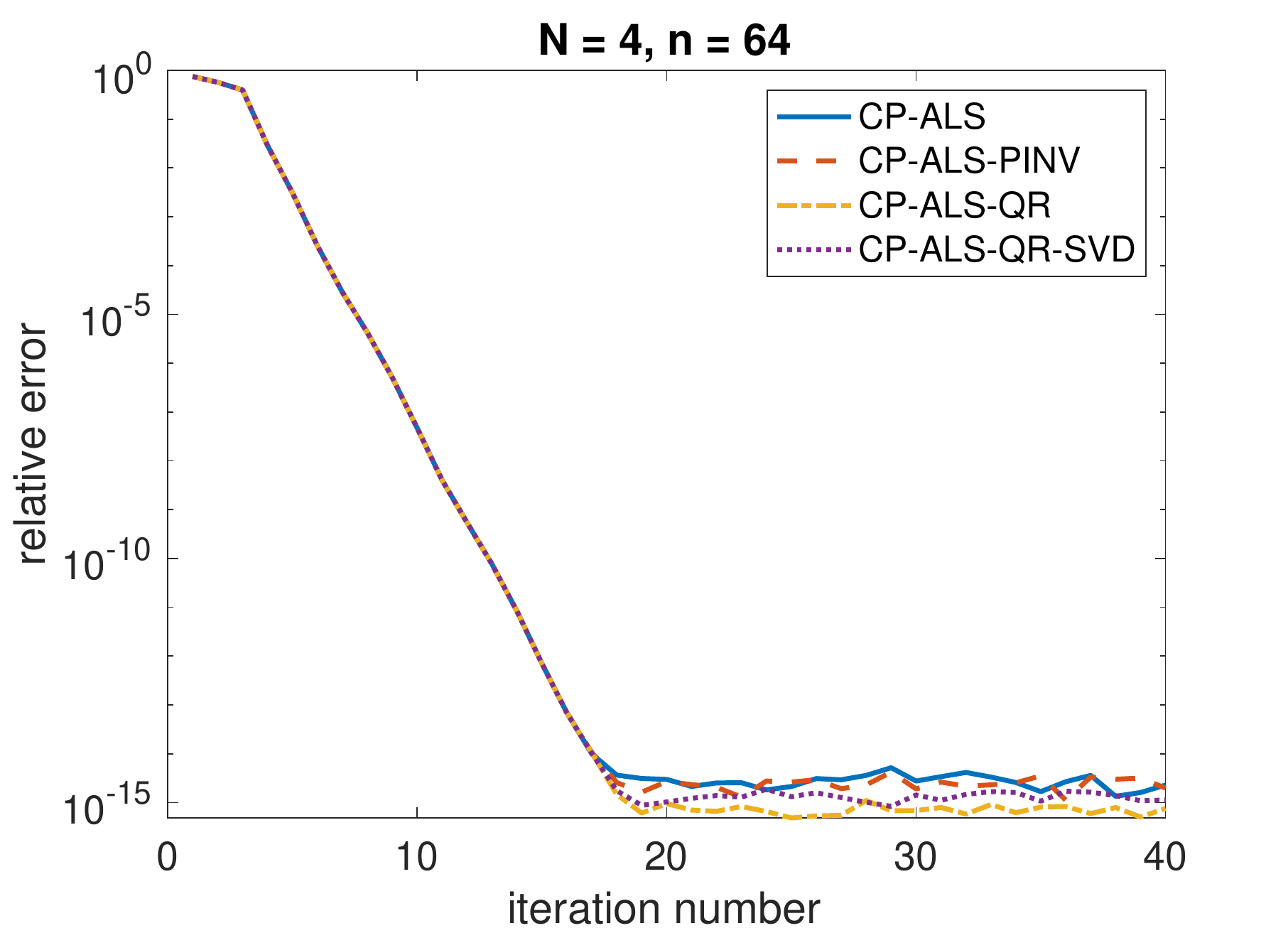}
\includegraphics[scale=.4]{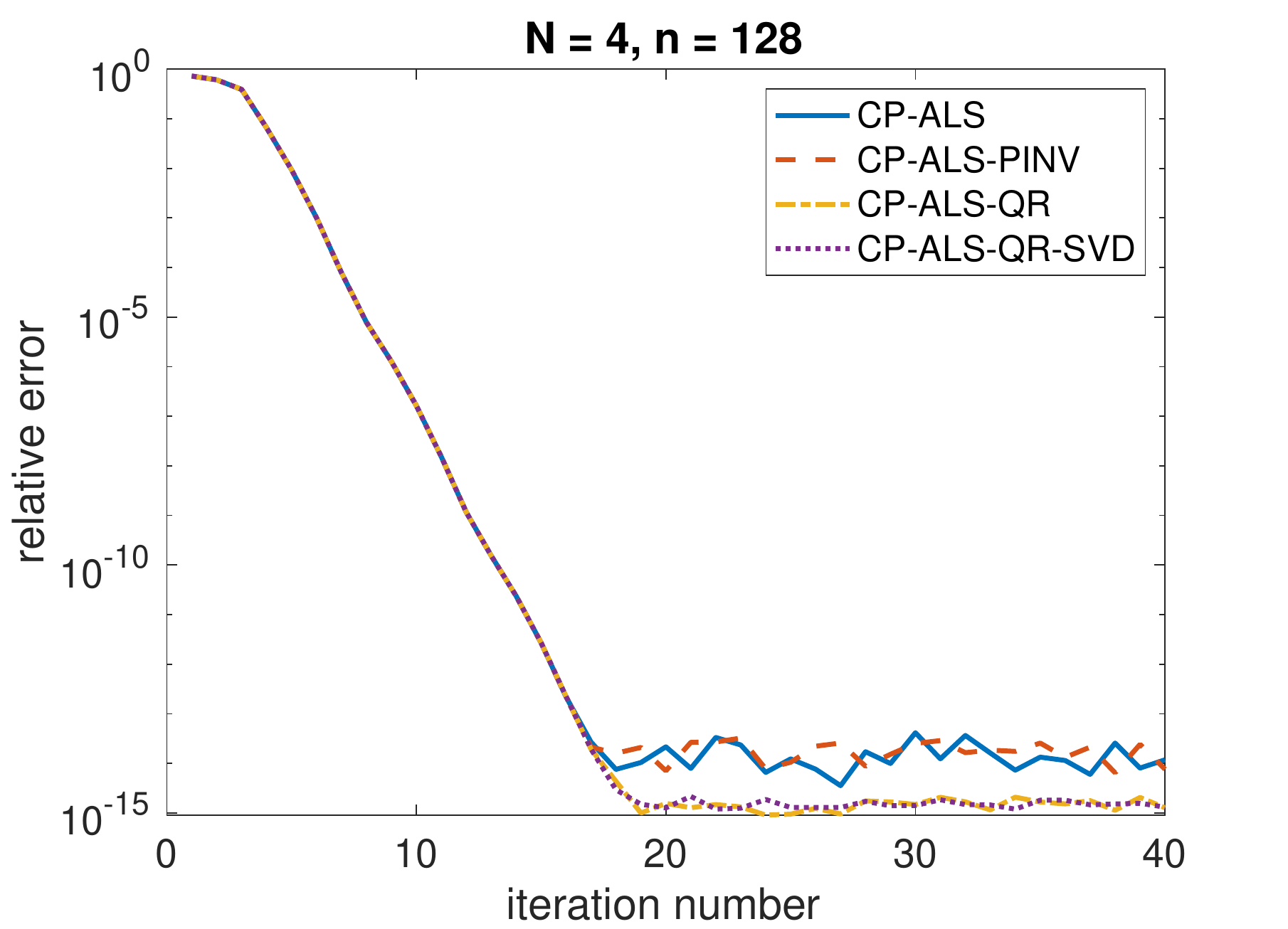}
\caption{Relative error of CP-ALS, CP-ALS-PINV, CP-ALS-QR, and CP-ALS-QR-SVD on the four-way sine of sums tensor in full rank-8 representation at the end of each iteration for dimension $n=64$ (left) and $n=128$ (right).}
\label{fig:sin4d}
\end{figure}

We first examine two lower-order cases, with $N = 4$ and 5 modes. 
Starting with $N = 4$, we compute the relative error at the end of each iteration for all four ALS algorithms for two different $n$ values, the size of vectors $x_j$. The results are plotted in \cref{fig:sin4d}, where we see that as $n$ increases, the difference between the relative errors of the algorithms becomes larger. Specifically, the relative error at the end of each iteration is larger for CP-ALS and CP-ALS-PINV than our QR-based algorithms. 

\begin{figure}
\centering
\includegraphics[scale=.4]{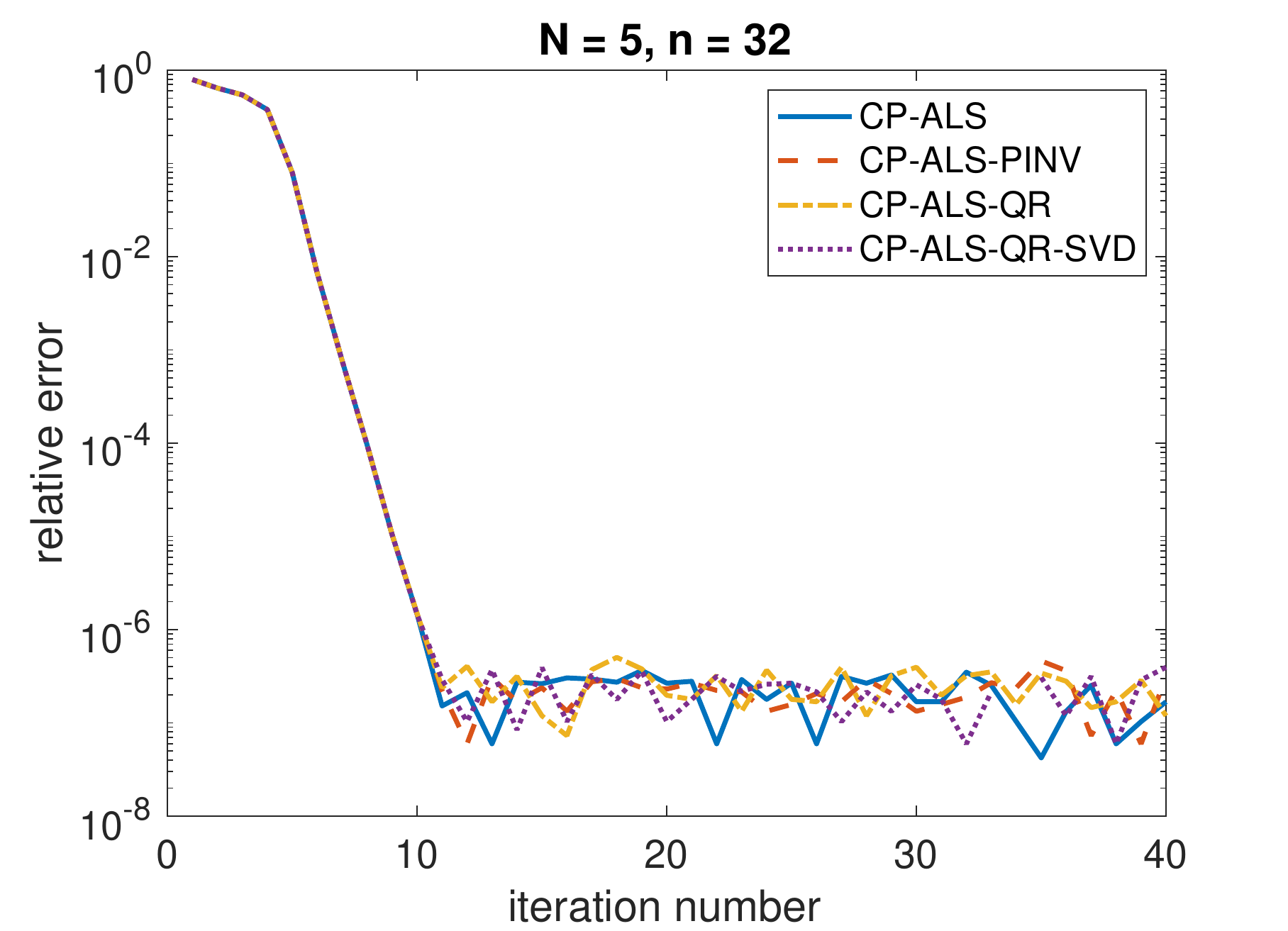}
\includegraphics[scale=.4]{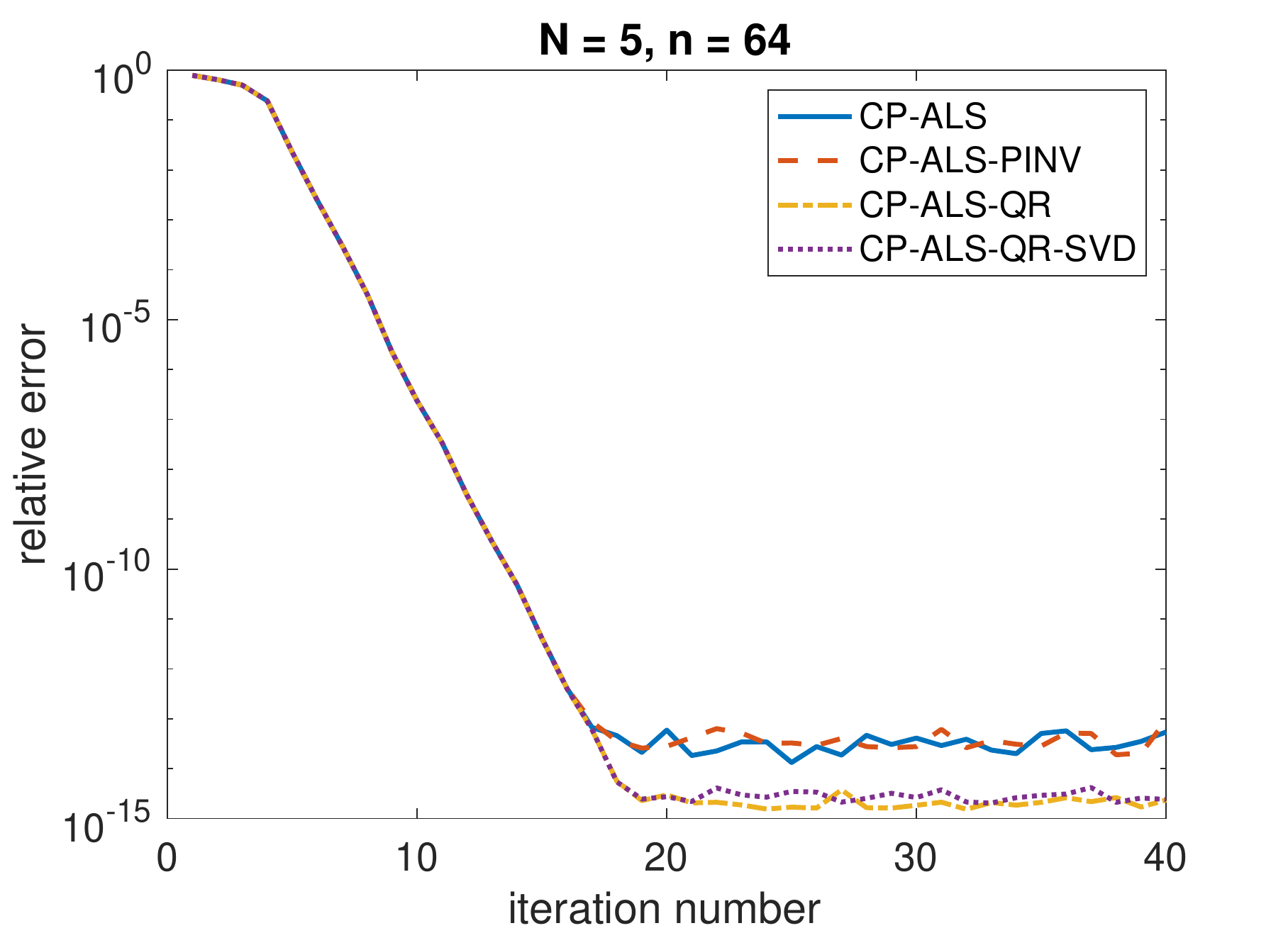}
\caption{Relative error of CP-ALS, CP-ALS-PINV, CP-ALS-QR, and CP-ALS-QR-SVD on the five-way sine of sums tensor in full rank-16 representation at the end of each iteration for dimension $n=32$ (left) and $n=64$ (right).}
\label{fig:sin5d}
\end{figure}

Similarly, for $N = 5$, we plot the relative error at the end of each iteration for all four ALS algorithms for increasing $n$ values. 
These results are shown in \cref{fig:sin5d}, where we see a similar result to $N = 4$. 
As the dimension gets bigger, the QR-based algorithms converge to a lower relative error value.
Note that in this case, $n=32$ is not a fine enough discretization to obtain relative error near machine precision; only $10^{-6}$ error is achieved by any algorithm.

\subsubsection{Higher-order tensors} 
\begin{figure}
\centering
\includegraphics[scale=.4]{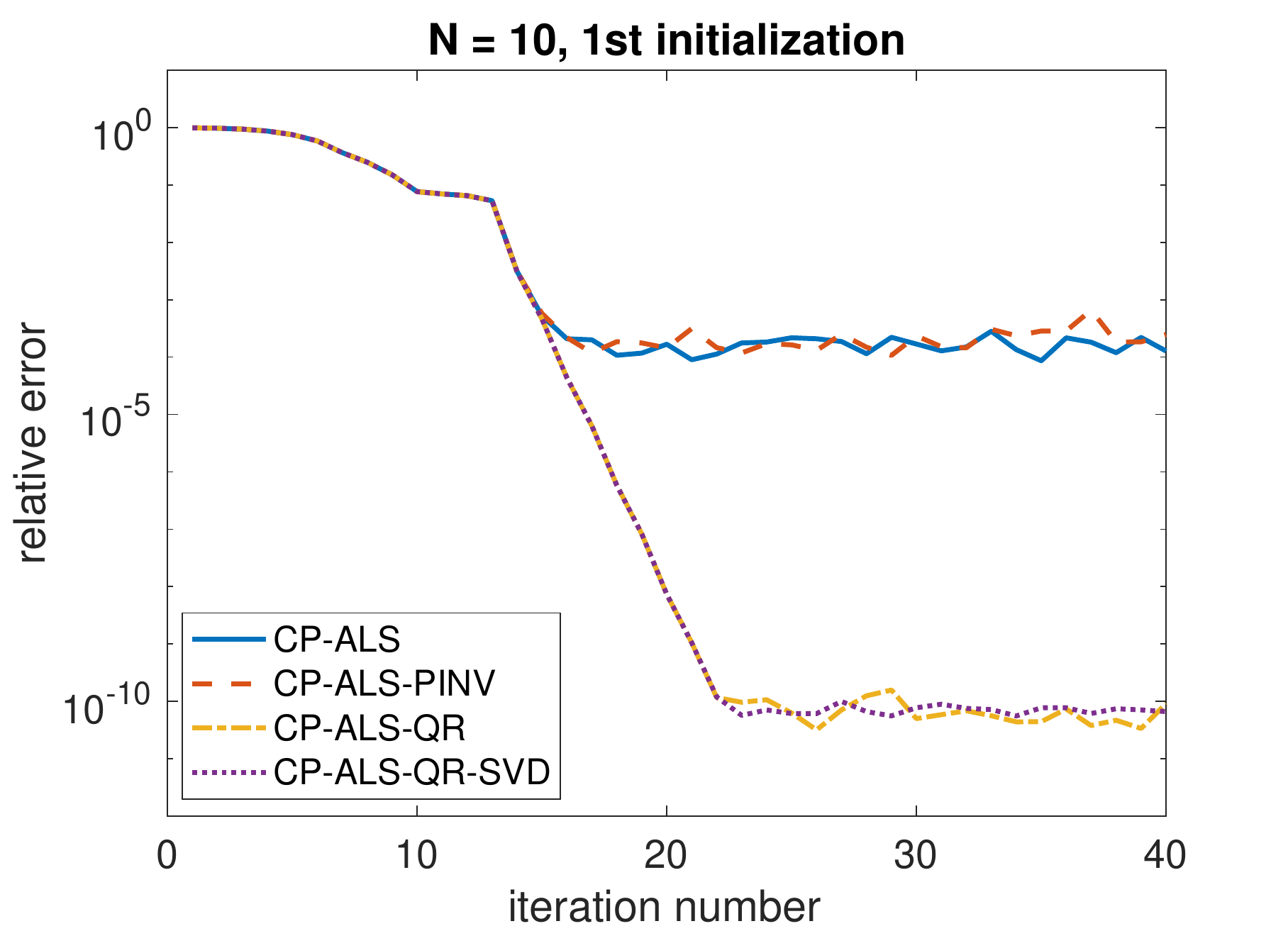}
\includegraphics[scale=.4]{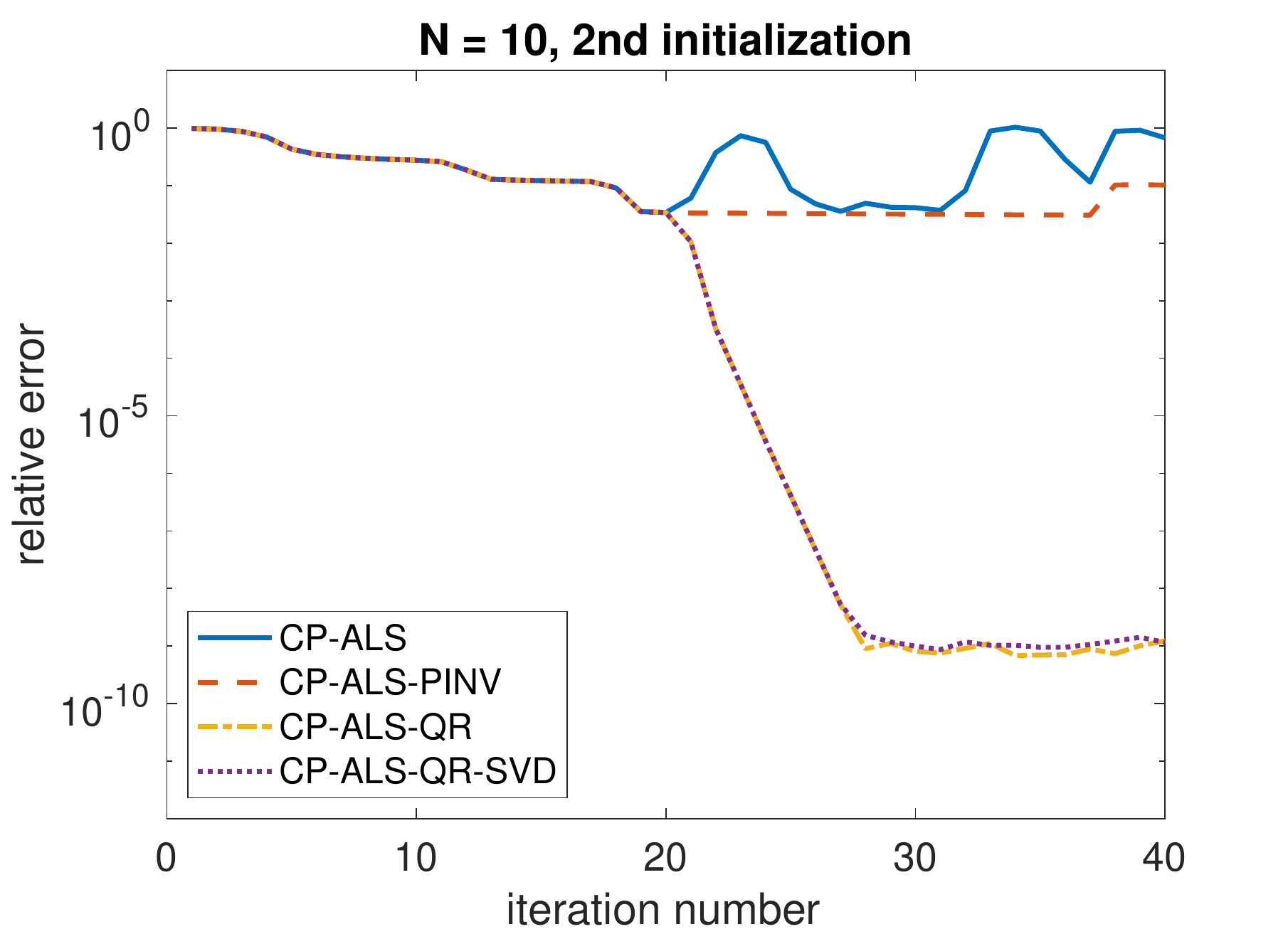}
\caption{Relative error of CP-ALS, CP-ALS-PINV, CP-ALS-QR, and CP-ALS-QR-SVD on the ten-way sine of sums tensor in rank-10 representation at the end of each iteration for dimension $n=8$. Results from two different random initializations are shown. }
\label{fig:sin10d}
\end{figure}

We now consider a higher-order case, where our tensor has $N = 10$ modes. 
In \cref{fig:sin10d}, we plot these relative errors in two cases. We use $n=8$ for both cases, but use two different random initializations to show two different types of results we obtained. For the first random initialization (left), we see similar results to the lower-order cases, with all four algorithms converging, but the gap between the relative error for the QR-based algorithms and the normal equations-based algorithms is much larger than in lower-order tensors. With the second random initialization (right), CP-ALS and CP-ALS-PINV do not converge to anything, while the relative error for CP-ALS-QR and CP-ALS-QR-SVD converge normally to low values. When repeating this experiment for multiple random initializations, we found that this second case was more common, occurring in four out of five trials. 

From these experiments, we can see that for ill-conditioned problems, our QR-based algorithms are more stable than the typical algorithms in higher dimensions. For all dimensions, we are also able to attain higher accuracy than traditional ALS algorithms.

\section{Conclusions and Future Work}
\label{sec:conclusion} 

We have developed and implemented versions of the CP-ALS algorithm using the QR decomposition and SVD in an effort to address the numerical ill-conditioning to which the normal equations in the traditional algorithm are susceptible. 
The first version uses a QR factorization to solve the linear least squares problems within CP-ALS. 
We also present the CP-ALS-QR-SVD algorithm, which applies the SVD as an extra step in the algorithm to handle numerically rank-deficient problems. 
In addition to the algorithms themselves, we provide analysis of their complexity, which is comparable to that of the widely-used CP-ALS algorithm when the rank is small. 
Our new algorithms prove useful for computing CP tensor decompositions with more stability in the event of ill-conditioned coefficient matrices, and present an alternative when analyzing tensor data for which the CP-ALS algorithm produces dissatisfactory results, or is unable to produce any result due to ill-conditioning. 
We envision our QR-based algorithms being used as part of a robust CP-ALS solver that uses the traditional normal equations approach by default, but solves the least squares problems using QR if any ill-conditioning is detected.

There are several potential performance improvements to pursue in future work. 
In situations where the target rank is high, computing $\M{Q}_0$, which involves a QR of a Khatri-Rao product of upper triangular matrices, becomes a more dominant cost of the CP-ALS-QR algorithm. 
The Khatri-Rao product of upper triangular matrices has structure which we do not exploit in our implementation, and which could lead to a more efficient implementation.  
We could also use dimension trees to speed up our implementation of the Multi-TTM function, the major dominant cost in both new algorithms. 
We currently use the Tensor Toolbox implementation but could improve the performance by reusing computations as in a dimension tree. 
The Gauss-Newton algorithm we use solves the approximate linear least squares problems via the normal equations. 
Another interesting direction to pursue would be to use the QR of the Jacobian to solve the least squares problems instead to improve the stability in the presence of ill-conditioning.


\bibliographystyle{siamplain}
\bibliography{main}

\end{document}

%% file: tensor_header.tex
\usepackage{amssymb,amsfonts,amsmath}

\usepackage[mathscr]{eucal}

\usepackage{stmaryrd}

\usepackage{amsbsy}

\usepackage{bm}

\usepackage{braket}

\usepackage{mathtools}

\newcommand{\V}[2][]{{\bm{#1\mathbf{\MakeLowercase{#2}}}}} 
 
\newcommand{\M}[2][]{{\bm{#1\mathbf{\MakeUppercase{#2}}}}} 
\newcommand{\Mn}[3][]{{\bm{#1\mathbf{\MakeUppercase{#2}}}}_{#3}} 
 
\newcommand{\T}[2][]{\boldsymbol{#1\mathscr{\MakeUppercase{#2}}}} 

\newcommand{\Mz}[3][]{\M[#1]{#2}_{(#3)}}

\newcommand{\Kron}{\otimes} 
\newcommand{\Khat}{\odot} 
\newcommand{\Hada}{\ast} 
 
\newcommand{\mb}[1]{\mathbb{#1}}
\newcommand{\bmat}[1]{\begin{bmatrix} #1 \end{bmatrix}}